\documentclass{article}
\usepackage{PRIMEarxiv}

\usepackage[utf8]{inputenc} 
\usepackage[T1]{fontenc}    
\usepackage{hyperref}       
\usepackage{url}            
\usepackage{booktabs}       
\usepackage{amsfonts}       
\usepackage{nicefrac}       
\usepackage{microtype}      
\usepackage{lipsum}
\usepackage{fancyhdr}       
\usepackage{graphicx}       
\graphicspath{{media/}}     

\usepackage{siunitx}    
\usepackage{comment}
\usepackage[export]{adjustbox}
\usepackage{hyperref}
\usepackage{makecell}
\usepackage{url}
\usepackage{tikz, pgfplots}
\usetikzlibrary{positioning}
\usepackage{capt-of}
\usepackage{diagbox}

\usepackage[ruled,vlined]{algorithm2e}
\usepackage{algpseudocode}
\def\eqref#1{(\ref{#1})}

\usepackage{amsfonts,bm}

\def\eqref#1{equation~\ref{#1}}

\def\1{\bm{1}}

\DeclareMathAlphabet{\mathsfit}{\encodingdefault}{\sfdefault}{m}{sl}
\SetMathAlphabet{\mathsfit}{bold}{\encodingdefault}{\sfdefault}{bx}{n}

\newcommand{\X}{{\mathbf X}}

\newcommand{\x}{{\boldsymbol x}}

\newcommand{\s}{{\boldsymbol s}}
\newcommand{\y}{{\boldsymbol y}}

\newcommand{\Y}{{\mathbf Y}}
\newcommand{\z}{{\boldsymbol z}}
\newcommand{\n}{{\boldsymbol n}}

\newcommand{\w}{{\boldsymbol w}}

\newcommand{\Rd}{{\mathbb R}}

\newcommand{\Ed}{{\mathbb E}}

\newcommand{\Ac}{{\mathcal A}}

\usepackage[utf8]{inputenc} 
\usepackage[T1]{fontenc}    
\usepackage{hyperref}       
\usepackage{url}            
\usepackage{booktabs}       
\usepackage{amsfonts}       
\usepackage{nicefrac}       
\usepackage{microtype}      
\usepackage{xcolor}       
\usepackage{xkcdcolors}
\usepackage{times}
\usepackage{epsfig}
\usepackage{graphicx}
\usepackage{placeins}
\usepackage{multirow}
\usepackage[skip=5pt]{caption}
\usepackage{rotating}
\usepackage{cancel}
\usepackage{setspace}
\usepackage{eso-pic}
\usepackage{amsmath}
\usepackage{bm}
\usepackage{bibunits} 

\usepackage{amssymb}
\usepackage{hyperref}
\usepackage{graphicx}
\usepackage{wrapfig,lipsum,booktabs}

\usepackage{epsfig}
\usepackage{tikz}
\usetikzlibrary{spy}
\usepackage{mathrsfs}

\usepackage{nicefrac}       
\usepackage{booktabs}       

\usepackage{thmtools,thm-restate}
\usepackage{cleveref}

\usepackage{subcaption}        

\title{Diff-ANO: Towards Fast High-Resolution Ultrasound Computed Tomography via Conditional Consistency Models and Adjoint Neural Operators
}

\author{
  Xiang Cao$^{a}$, Qiaoqiao Ding$^{a}$, Xinliang Liu$^{b}$, Lei Zhang$^{a}$,  Xiaoqun Zhang$^{a}$ \\
  $^{a}$ School of Mathematical Sciences and Institute of Natural Sciences \\
  Shanghai Jiao Tong University \\
  Shanghai, 200240, CHINA \\ 
  $^{b}$ School of Mathematical Sciences, Ocean University of China \\
  Qingdao, 266100, CHINA \\
}

\begin{document}
\maketitle

\begin{abstract}
Ultrasound Computed Tomography (USCT) involves reconstructing medium properties from scattered acoustic waves. At its core, this is a PDE-constrained optimization problem minimizing a data-fidelity objective, whose inherent ill-posedness severely challenges traditional numerical solvers. While standard approaches struggle with skip-cycle local minima and computational bottlenecks, generative models applied to PDE inverse problems offer a powerful data-driven initialization and regularization strategy. In this work, we introduce \textbf{Diff-ANO} (\textbf{Diff}usion-based Models with \textbf{A}djoint \textbf{N}eural \textbf{O}perators), a framework shifting the focus toward generative models applied to PDE inverse problems by combining consistency modeling with adjoint operator learning. Rather than performing a fully Bayesian posterior sampling, our method efficiently solves the measurement-constrained least-squares problem through two key innovations: (1) a \textit{conditional consistency model} that provides a high-quality data-informed prior to initialize the optimization, mapping directly from noisy states to the clean manifold, and (2) an \textit{adjoint operator learning} module that replaces slow PDE solves with neural operator surrogates for rapid measurement-informed gradient computation. To enable this, we introduce the batch-based Convergent Born Series (BCBS)—accelerating the production of wavefield data for training the neural adjoint operator. Experiments demonstrate significant improvements in both efficiency and reconstruction quality, particularly under sparse-view and partial-view scenarios.
\end{abstract}


\section{Introduction}\label{sec:1}

\subsection{Ultrasound Computed Tomography (USCT)}\label{sec:1.1}
USCT formulates a nonlinear, PDE-constrained inverse problem that aims to recover the sound-speed distribution within a medium from the acoustic wave measurements. This modality has seen widespread application in medical imaging and geophysical exploration for high-resolution tomographic images \cite{arridge2009optical}. Mathematically, the forward model for USCT is governed by the Helmholtz equations: for a fixed  angular frequency $\omega$, the acoustic wavefield $\Y_n(\mathbf{r})$ for each source point in $\{\mathbf{r}^s_n\}_{n = 1}^{N}$ satisfies
\begin{equation}\label{helmholtz}
    \nabla^2 \Y_n(\mathbf{r}) + \frac{\omega^2}{\X_0(\mathbf{r})^2} \Y_n(\mathbf{r}) = -\rho_n(\mathbf{r}), \; \forall \mathbf{r}\in \Omega, 
\end{equation} 
where $\X_0(\mathbf{r})$ denotes the sound-speed distribution and $\rho_n(\mathbf{r}) := \delta(\mathbf{r}-\mathbf{r}^s_n)$ is the Dirac delta-function. Typically, $\X_0$ is non-uniform only within $\Omega_0$, a predefined domain of interest (DOI), while the surrounding region $\Omega$\textbackslash$\Omega_0$ has constant background speed. Absorbing boundary conditions are widely used to emulate wave absorption at the medium boundaries \cite{colton2013inverse}. 

\begin{figure}
    \centering
    \includegraphics[width=0.875\linewidth]{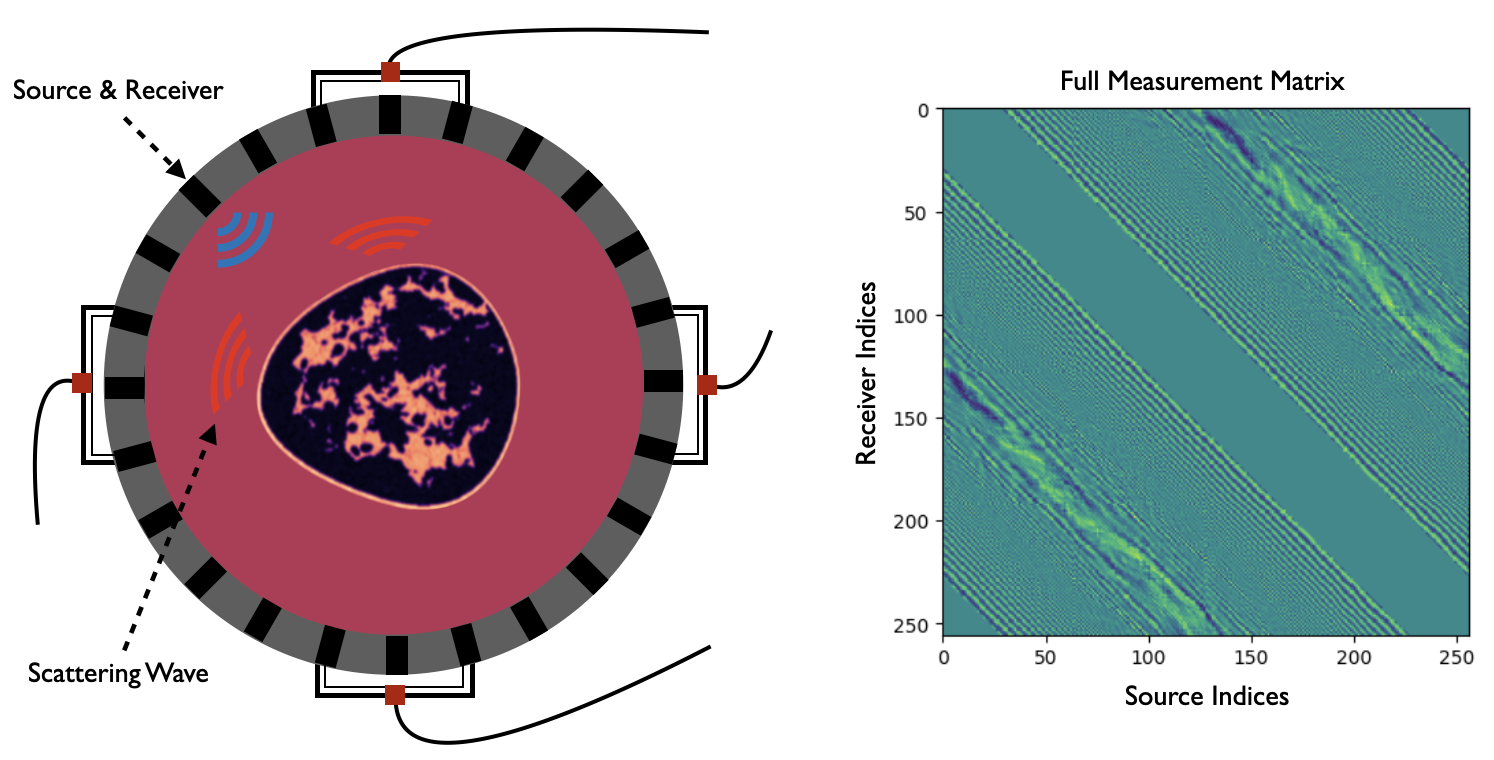}
    \caption{Schematic of the USCT measurement setup and example data. \emph{Left:} Circular array of transmitters and receivers surrounding the region of interest, illustrating wave emission, scattering, and reception. \emph{Right:} Real part of the full measurement data matrix $\y$, with source indices on the horizontal axis and receiver indices on the vertical axis.}
    \label{fig:usct-physics-illustration}
\end{figure}

Denote by $\mathcal{K}:(\X_0 ;\rho_n)\in \mathcal{X}(\Omega) \times \mathcal{X}(\Omega) \rightarrow \Y_n \in \mathcal{Y}(\Omega)$ the Helmholtz solution operator defined by \cref{helmholtz}, where $\mathcal{X}(\Omega)$ and $\mathcal{Y}(\Omega)$ are two Banach function spaces. In what follows, we assume receivers and sources are co-located, i.e. $\{\mathbf{r}_m\}_{m = 1}^{M} = \{\mathbf{r}^s_n\}_{n = 1}^{N}$. The associated inverse problem seeks to estimate $\X_0$ from noisy measurements $\y^\delta \in \mathbb{C}^{M \times N}$, with its $(m,n)$-th entry defined as 
\begin{equation}
    \y_{m,n}^\delta =  \mathcal{K}(\X_0;\rho_n)(\mathbf{r}_m) + \delta_{m,n} \eta,
\end{equation}
where $\eta$ is a complex random noise, and $\delta_{m,n}$ controls the noise level for each source-receiver pair. This inverse problem is severely ill-posed due to three principal factors:

\begin{enumerate}
  \item \textit{Measurement Noise.} Sensor-model discrepancies and environmental noise corrupt each receiver’s measurement, degrading resolution and biasing reconstructions unless properly modeled. 
  \item \textit{Incomplete Data Scenarios.} Analogous to sparse‐view and partial‐view scenarios in X-ray Computed Tomography, limited angular coverage or receiver count in USCT severely degrades the conditioning of the inversion \cite{openpros2025}.
  \item \textit{Skip-Cycle Phenomena.}  Multiple scattering and resonance can introduce cycle
  skipping in phase‑based inversion, leading to convergence to false local minima and destabilizing the inversion if not properly addressed \cite{pratt1998gauss}. 
\end{enumerate}
The schematic illustration is shown in \cref{fig:usct-physics-illustration}. To mitigate these challenges, one must incorporate robust priors or regularization—such as variational Bayesian formulations~\cite{fortuin2021priors}, total-variation and sparsity priors \cite{li2018adaptive} to exploit spatial similarity. Other advanced learning-based approaches in USCT including plug-and-play priors \cite{kamilov2022pnp}, untrained neural representations \cite{yan2024plug} and generative diffusion priors \cite{wang2023prior} to ensure stable and accurate sound-speed reconstructions. 

\subsection{Solving Inverse Problems Using Diffusion Models}\label{sec:1.2}
The inherent ill-posedness of inverse problems necessitates a dual emphasis on \emph{data fidelity} and \emph{prior regularization} to stabilize solutions \cite{arridge2019solving}. In the context of USCT, one fundamentally aims to solve a least-squares optimization problem to find a sound-speed model consistent with the measured data. However, the severe nonlinearity and limited data make this optimization highly non-convex.

To regularize this process, generative models have emerged as powerful tools. Rather than viewing the inversion strictly through a fully Bayesian lens where one exhaustively samples the posterior distribution, we treat the generative model as a sophisticated, data-informed prior. It serves to provide a robust initialization and regularization for the least-squares problem, projecting intermediate solutions into a manifold of physically realistic sound-speed maps. Within a statistical framework \cite{tarantola2005inverse}, the prior $p(\x_0)$ guides the likelihood $p(\y^\delta|\x_0)$:
\begin{equation}
    p(\x_0|\y^\delta) \propto \underbrace{p(\y^\delta|\x_0)}_{\text{Likelihood}} \cdot \underbrace{p(\x_0)}_{\text{Prior}}, 
\end{equation}
where the prior $p(\x_0)$ regularizes solutions by encoding domain-specific knowledge—a component historically limited by handcrafted designs (e.g., sparsity or total variation \cite{rudin1992nonlinear}). 

Recent advances in diffusion models \cite{ho2020denoising, yang2023diffusion} have revolutionized this paradigm by learning \emph{implicit data-driven priors} through iterative denoising processes. These models approximate the unconditional score function $\nabla_{\x_t} \log p_t(\x_t)$, which guides sampling trajectories toward high-probability regions of the data manifold. For inverse problems, diffusion-based posterior sampling leverages a conditional reverse process derived via Bayes’ rule:
\begin{equation}
\nabla_{\x_t} \log p_t(\x_t|\y^\delta) = \underbrace{\nabla_{\x_t} \log p_t(\y^\delta|\x_t)}_{\text{Likelihood gradient}} +  \underbrace{\nabla_{\x_t} \log p_t(\x_t)}_{\text{Prior score}},
\end{equation}
where the prior score is learned from data, while the likelihood gradient ties measurements to the latent variable $\x_t$.

For applications, most research primarily focuses on linear/nonlinear inverse problems in non-PDE contexts, where the forward operator $\Ac$ can be explicitly formulated as differentiable compositions amenable to automatic differentiation. For linear inverse problems—such as inpainting \cite{lugmayr2022repaint}, deblurring \cite{chung2023diffusion}, super-resolution \cite{li2022srdiff}, Computed Tomography \cite{song2022solving}, and Magnetic Resonance Imaging \cite{chung2022score}—the forward model often reduces to structured matrix operations (e.g., Radon transforms in CT). Similarly, nonlinear problems like phase retrieval \cite{peer2023diffphase}, nonlinear deblurring \cite{chung2023diffusion}, and high dynamic range imaging \cite{mardani2023variational} leverage differentiable physics-inspired models. These frameworks benefit from gradient-based optimization, where the likelihood term $\nabla_{\x_t} \log p_t(\y^\delta|\x_t)$ is efficiently approximated via backpropagation through $\mathcal{A}$. However, extending these methodologies to nonlinear PDE-constrained inverse problems—such as USCT governed by the Helmholtz equation—faces three principal challenges:

\begin{enumerate}
    \item \textit{PDE-Constrained Gradient.} The inherent nonlinearity of the forward operator $\Ac$ introduces a critical dependency of its Fréchet derivative $(\partial \Ac)_{\x_0}$ on the parameter field $\x_0$. This necessitates real-time computation of the Jacobian-vector product (JVP) to evaluate the data fidelity gradient: 
    \begin{equation*}\label{constraints}
        \nabla_{\x_0} \left\|\y^{\delta} - \mathcal{A}(\x_0)\right\|_2^2 = 2(\partial \Ac)^{*}_{\x_0}\left(\mathcal{A}(\x_0) - \y^{\delta}\right)
    \end{equation*} 
    where $(\partial \Ac)^{*}_{\x_0}$ denotes the adjoint operator. While automatic differentiation efficiently computes gradients for explicit forward models, it falters in PDE-based systems due to the \textit{implicit} coupling between $\mathcal{A}(\x_0)$ and $\x_0$. For example, solving the Helmholtz equation iteratively embeds $\x_0$ into the solver’s internal states, precluding direct AD-based differentiation. 

   \item \textit{Discretization-Induced Approximation Error.} In PDE-based inverse problems, the governing PDEs (e.g., Helmholtz equations) are inherently formulated in the continuous domain, while diffusion sampling operates on discretized grids. Discretization of the PDE solution operator $\mathcal{K}$—via finite element methods or finite differences—introduces numerical approximation errors that propagate through multi-step sampling. It is necessary to theoretically bridge the gap from misaligned domains between continuous PDE formulations for $\X_0$ and discrete score networks for $\x_0$. 
    
   \item \textit{Computational Imbalance.} Diffusion-based posterior sampling requires hundreds to thousands of sequential evaluations of the score function $\nabla_{\x_t} \log p_t(\x_t)$ and data fidelity gradients. Each evaluation step demands solving the related forward or adjoint PDEs under multiple boundary conditions, drastically slowing inference compared to unconditional sampling \cite{cao2025subspace}. The computational complexity will scales with the dimensionality of the PDE discretization, and the iterative nature of numerical solvers (e.g., Convergent Born Series for Helmholtz \cite{osnabrugge_convergent_2016}). While score models benefit from GPU acceleration, well-established CPU-based PDE solvers  will dominate runtime, creating a hardware mismatch. 
\end{enumerate}

\subsection{Main Work} \label{sec:1.3}

\begin{figure*}[!htp]
    \centering
    \includegraphics[width=\linewidth]{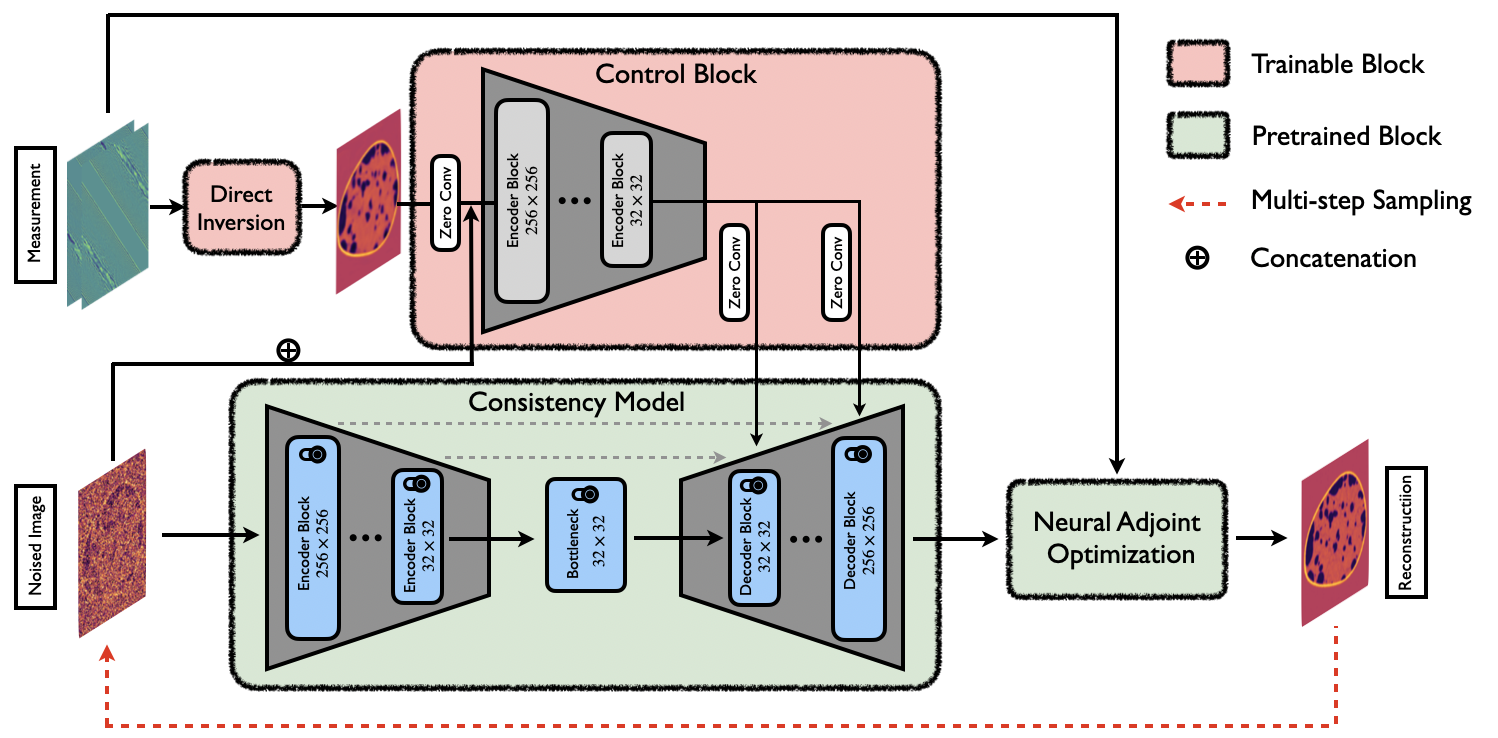}
    \caption{Proposed pipeline of the \textbf{Conditional Consistency Model with Adjoint Neural Operator} for USCT. Measurements are first mapped to a coarse image by the \emph{Direct Inversion} module (pink). This provisional reconstruction drives a trainable \emph{Control Block} (pink) that modulates a pretrained \emph{Consistency Model} (green) via zero-convolution adapters. The consistency model performs multi-step sampling (red dashed path) to progressively refine the image. Meanwhile, the \emph{Neural Adjoint Optimization} module  (green) enforces physics constraints via a neural operator surrogate of the Helmholtz operator.}
    \label{fig:framework}
\end{figure*}

Our motivation lies in replacing the traditional adjoint-based optimization, which typically relies on computationally intensive numerical PDE solvers, with a more efficient \textit{neural operator} approach in the consistency sampling framework.  First, the ill-posed nature of USCT demands an appropriate initial reconstruction for optimization. Here, the conditional \textit{consistency model} , as a generalization of the direct supervised model, provides iterative refinement that better aligns the reconstruction with the data prior. This alignment is achieved by conditioning the sampling process of the consistency model on the direct inversion as \cite{zhao2024cosign}. Second, PDE-based inverse problems typically require adjoint-based optimization using well-established numerical solvers. By exploiting the self-adjointness of the Helmholtz operators, we integrate the pretrained neural operator into the adjoint-based optimization, ensuring measurement consistency across the multi-step sampling process.

To enable practical deployment, we introduce the batch-based \textit{Convergent Born Series}~(BCBS), a memory-efficient strategy for on-the-fly generation of neural operator training pairs. Notably, our neural operators require training only on the clean data manifold rather than the full optimization trajectory, substantially improving generalization while reducing sample complexity.  Comprehensive numerical experiments demonstrate that our designed framework achieves rapid and high-fidelity USCT reconstructions in few-step evaluations by simultaneously enforcing: (1) physics constraints through neural operator surrogates, and (2) data priors through the conditional consistency model. Here, the unified architecture is illustrated in \cref{fig:framework}. 

\paragraph{Organization} The remainder of this paper is organized as follows. In \cref{sec:2}, we review related works on diffusion models for inverse problems and neural operators for PDE solvers. \cref{sec:3} presents the fundamentals of score-based diffusion models and their distillation into consistency models, and also introduces Diffusion Posterior Sampling for solving inverse problems. Building upon these foundations, \cref{sec:4} details our proposed conditional consistency model with adjoint neural operators for USCT, where BCBS is adopted for online training. In \cref{sec:5}, we describe implementation specifics, including measurement configurations and network training settings. \cref{sec:6} compares our method against baselines and reports numerical results. Furthermore, we conduct the ablation study from two aspects—inversion blocks and forward neural operators—to validate the designed components in \cref{sec:7}. Finally, \cref{sec:8} discusses the limitations of our approach and concludes the paper.

\section{Related work}\label{sec:2}
Our work primarily focuses on solving inverse problems using diffusion models and accelerating solving PDEs via neural operators. We will briefly review related works about these topics in the following content. 

\textbf{Diffusion Models for Inverse Problems.} 
The first category leverages Bayesian inference frameworks to estimate posterior distributions conditioned on measurements. A foundational approach uses the Tweedie formula \cite{robbins1992empirical} or measurement subspace projections \cite{kadkhodaie2021stochastic} to guide conditional sampling. Building on this, {Pseudoinverse-guided Diffusion Models} ($\Pi$GDM) \cite{song2023pseudoinverse} introduces pseudo-inverse guidance for linear inverse problems, achieving exact consistency for tasks like image inpainting. To enable zero-shot restoration without retraining, {Denoising Diffusion Nullspace Model} (DDNM) \cite{wang2022zero} decomposes the measurement-consistent solution into range-space and null-space components. While effective for linear forward cases, these methods struggle with nonlinear inverse problems. DPS \cite{chung2023diffusion} addresses this by incorporating likelihood gradients via automatic differentiation, enabling applications to general nonlinear operators. Decomposed Diffusion Sampler (DDS) \cite{chung2023decomposed} enables conjugate gradient  optimization on Tweedie-denoised samples, eliminating the need for manifold-constrained gradient~(MCG) computations \cite{chung2022improving}. In DDS, the DDIM \cite{song2020denoising} sampling acceleration can be applied to further expedite the posterior sampling process of DDS. 

The second category reformulates inverse problems as variational optimization with diffusion priors. {RED-diff} \cite{mardani2023variational} bridges denoising-based regularization \cite{romano2017little} and diffusion models, casting sampling as stochastic optimization of a RED-inspired loss. From a Bayesian filtering perspective, FPS \cite{dou2024diffusion} establishes theoretical guarantees for diffusion-based inverse problem solving by revealing the equivalence between posterior sampling and sequential Monte Carlo filtering. These methods are problem-agnostic, requiring no task-specific training. In contrast, problem-specific approaches train conditional diffusion models for targeted applications. For instance, \cite{whang2022deblurring} trains a deblurring-specific diffusion model using paired datasets, while \cite{rissanen2022generative} learns the inverse heat dissipation process as a diffusion model for heat-inversion. 

A comprehensive taxonomy of these methods is provided in \cite{daras2024survey}, highlighting their applicability to non-PDE inverse problems. However, diffusion-based approaches for PDE-constrained inverse problems remain underexplored due to challenges in enforcing the PDE-constrained gradient and handling the computational imbalance, as discussed in \cref{sec:1.2}

\textbf{Diffusion Sampling Acceleration.} 
Traditional diffusion models suffer from slow sampling due to hundreds of sequential steps in reverse SDEs. Many studies \cite{bao2022analytic, dhariwal2021diffusion} have focused on reducing the number of discretized sampling steps with adaptive solvers for the reverse process. Diffusion Probabilistic Model ODE-Solver (DPM-Solver) \cite{lu2022dpm}, as a generalization of DDIM \cite{song2020denoising}, solves the probability flow ODE using high-order numerical solvers. Despite these, model-based acceleration approaches mitigate the slow-sampling issue via architectural innovations: Subspace Diffusion Generative Models (SDGM) \cite{jing2022subspace} reduce computational costs for score evaluations by restricting the diffusion process through projections onto lower-dimensional subspaces. Rather than operating in pixel spaces, Latent Diffusion Models (LDMs) \cite{vahdat2021score} is designed in low-dimensional latent spaces to reduce computational complexity. Through distillation from pretrained diffusion models, CMs \cite{song2023consistency} achieve one-step generation by learning self-consistency mappings across diffusion trajectories. Among them, CMs’ multi-step refinement capability is particularly promising for inverse problems. CoSIGN \cite{zhao2024cosign} introduces conditional CMs with ControlNet \cite{zhang2023adding} guidance, enabling few-step reconstruction with hard measurement constraints. This aligns with our approach: by combining CMs’ fast sampling with neural operators for surrogate measurement constraints, we achieve efficient and high-quality USCT reconstruction. 

\textbf{Neural Operators for Solving PDEs.}
Traditional neural networks are designed to map between finite-dimensional Euclidean spaces, whereas operator learning aims to approximate mappings between infinite-dimensional function spaces governed by PDEs \cite{boulle2024mathematical}. Neural operators have emerged as a powerful paradigm to directly learn the solution operator of PDEs, achieving orders-of-magnitude acceleration compared to classical numerical solvers \cite{li2020fourier}. Two seminal architectures exemplify this concept:
(1) DeepONet \cite{lu2021deeponet}, which employs a branch-trunk architecture theoretically grounded in universal approximation theory for operators; (2) Fourier Neural Operator (FNO) \cite{li2020fourier}, which parameterizes integral kernels in Fourier space to efficiently capture global spectral patterns. While spectral‐type operators such as FNO excel at capturing global structures, they often struggle with local details—e.g., boundary information or high-frequency features. By contrast, architectures like UNets \cite{unet} naturally accommodate complex boundaries but suffer from parameter inefficiency and limited long-range dependency modeling. To bridge this gap, Liu \emph{et al.} proposed the Hierarchical Attention Neural Operator (HANO), which mitigates spectral bias by adaptively coupling information across scales via attention, thereby boosting accuracy on challenging multiscale benchmarks \cite{liu2022ht}. Multigrid-inspired neural operators (MgNet, MgNO) \cite{he2019mgnet, he2023mgno} address these limitations by combining multigrid principles with neural networks.  MgNet, MgNO and its adaptations have been explored for a broader class of numerical PDEs \cite{chen2022meta, he2023mgno}. The inherited multigrid structure ensures alignment with PDE discretization hierarchies, making this  multigrid-inspired backbone particularly suitable for wave-based PDEs, where both local scattering phenomenon and global wave propagation must be resolved. 

For inverse problems, neural operators have emerged as powerful tools for PDE-governed inverse problems, primarily through two paradigms: direct data-to-parameter mapping \cite{molinaro2023neural, guotransformer, liu2025neumann} and accelerated forward/adjoint modeling for iterative optimization. While Neural Inverse Operators (NIOs) \cite{molinaro2023neural} combining DeepONets and FNOs demonstrate impressive reconstruction speed, their performance is fundamentally constrained by ill-posedness arising from sparse or noisy measurements. To address this limitation, recent works adopt neural operators as surrogates for forward/adjoint PDE solvers, enabling efficient gradient-based inversion. This paradigm has achieved notable success in seismic inversion \cite{yang2023rapid} and USCT \cite{zeng2023neural, zeng2025openbreastusbenchmarkingneuraloperators, wu2026diffsos}. Crucially, inverse problems impose stricter requirements on neural operators compared to forward modeling. First, the trained operator must maintain high accuracy not merely on clean data distributions but across the entire optimization trajectories. Second, the adjoint operators’ structural dependence on forward solutions demands co-designed neural approximations. These aspects reveal that, the synergistic integration of two critical components—neural operators for accelerating forward/adjoint Helmholtz solves and diffusion-based priors for mitigating USCT's ill-posedness—remains an open challenge. 

\section{Preliminaries on Diffusion-Based Inverse Modeling}\label{sec:3}

\subsection{Score-Based Diffusion Models (SDMs) to Consistency Models  (CMs)
}\label{sec:3.1}
SDMs \cite{song2020score} characterize the forward noising of data $\x_t\in\mathbb{R}^d$ over $t\in[0,T]$ via the variance-preserving SDE
\begin{equation}\label{eq:forward-sde}
    d\x_t = -\frac{\beta(t)}{2}\x_t dt + \sqrt{\beta(t)}d\w,
\end{equation}
where $\beta(t)=\beta_{\min}+t(\beta_{\max}-\beta_{\min})$ schedules the diffusion strength and $\w_t$ is a standard $d$-dimensional Wiener process. As $t\to T$, $\x_t$ approaches $\mathcal{N}(\mathbf{0},\mathbf{I})$, providing a tractable terminal distribution. To sample from the data distribution $p_{\mathrm{data}}(\x_0)$, one reverses this process by solving the reverse-time SDE
\begin{equation}\label{eq:reverse-sde}
    d\x_t = \left[-\frac{\beta(t)}{2}\x_t - \beta(t)\nabla_{\x_t} \log p_t({\x_t})\right]dt + \sqrt{\beta(t)}d\bar\w,
\end{equation}
where $p_t({\x_t})$ denotes the marginal density and $\bar\w$ runs backward in time. The intractable score function $\nabla_{\x_t} \log p_t({\x_t})$ is approximated by a neural network $\s_{\Theta_{\mathrm{SD}}} (\x_t,t)$, trained via denoising score matching \cite{vincent2011connection}. 

Once the optimal $\Theta^*_{\mathrm{SD}}$ is learned, any numerical solver — e.g. Euler–Maruyama \cite{song2020denoising} or Predictor–Corrector \cite{nichol2021improved}, can plug in with $\s_{{\Theta^*_{\mathrm{SD}}}}(\x_t, t)$ to perform the reverse-time SDE ~\cref{eq:reverse-sde} and generate high-fidelity samples.

\textbf{Consistency Models}. Although SDMs achieve state-of-the-art sample quality, they typically require solving a multi-step reverse-time SDE/ODE, which can be computationally expensive. To address this, CMs \cite{song2023consistency} train a surrogate model for the score-based generative process. Specifically, the CM function $f_{\Theta_{\rm CM}}(\x_t, t)$ explicitly maps any intermediate noisy state $\x_t$ at time $t$ directly to the clean data manifold (the origin $\x_\varepsilon$, where $\varepsilon$ is a small positive constant).

\emph{Consistency Surrogate Modeling.} CMs can be obtained by surrogate modeling (often termed Consistency Distillation)—where a pretrained SDM's trajectory is learned by $f_{\Theta_{\rm CM}}$. This minimizes a self-consistency loss:
\begin{equation}\label{eq:consistency_loss}
    \mathcal{T}_{\rm CM}
     = d\bigl(f_{\Theta_{\rm CM}}(\x_{t_{n+1}}, t_{n+1}),\,f_{\Theta^{-}_{\rm CM}}(\widetilde{\x}_{t_n}, t_n)\bigr),
\end{equation}
where $\widetilde{\x}_{t_n}$ is obtained along the ODE trajectory, and $f_{\Theta^{-}_{\rm CM}}$ denotes the exponential moving average (EMA) as the target model.

\emph{Consistency Multi-Step Sampling.} While CMs excel at one-step generation, they also support an iterative multi-step sampler that refines sample quality. For a schedule $\varepsilon = \tau_0 < \tau_1 < \cdots < \tau_{S} = T$, one repeatedly applies:
\begin{equation}
    \begin{aligned}
        \x_0 &\leftarrow f_{\Theta^{*}_{\rm CM}}\bigl(\x_{\tau_{n+1}},\,\tau_{n+1}\bigr), \\
        \x_{\tau_n} &\leftarrow \text{Forward SDE}(\x_0,\tau_n), \quad n=S-1,\dots,0.
    \end{aligned}
\end{equation}
In words, this multi-step process alternates between two actions: first, the network projects the current noisy state onto the clean data manifold ($\x_0$); second, a controlled amount of noise is injected back into the clean sample to move it to the previous time step ($\x_{\tau_n}$). This loop progressively refines the sample, correcting approximation errors made in earlier single-step jumps.

\subsection{Diffusion Posterior Sampling for Inverse Problems}\label{sec:3.2} 
We consider the discrete forward model of an inverse problem in the general form:
\begin{equation}\label{eq:inverse}
  \y^{\delta} = \Ac(\x_0) + \bm{\n}, 
  \quad \y,\bm{\n}\in\Rd^n, \ \x_0\in\Rd^d,
\end{equation}
where \(\Ac(\cdot):\Rd^d\to\Rd^n\) denotes the forward operator and \(\bm{\n}\) is the additive measurement noise. Since recovering \(\x_0\) from \(\y^{\delta}\) suffers from the ill-posedness, we impose a learned diffusion prior and sample from the posterior via Bayes' theorem.  In continuous-time settings, the reverse–time SDE becomes
\begin{equation}\label{reverse-sde-posterior}
    \begin{aligned}
    d\x_t &= \left[-\frac{\beta(t)}{2}\x_t - \beta(t)(\nabla_{\x_t} \log p_t(\x_t|\y^{\delta}))\right]dt + \sqrt{\beta(t)}d\bar\w \\
    &= \left[-\frac{\beta(t)}{2}\x_t - \beta(t)(\nabla_{\x_t} \log p_t(\x_t) + \nabla_{\x_t} \log p_t(\y^{\delta}|\x_t))\right]dt + \sqrt{\beta(t)}d\bar\w, \\
    \end{aligned}
\end{equation} 
where the unconditional score function \(\nabla_{\x_t}\log p_t(\x_t)\) provides a plug-and-play prior for efficient posterior sampling. In DPS \cite{chung2023diffusion}, the intractable gradient of the data log-likelihood \(\nabla_{\x_t}\log p_t(\y^{\delta}|\x_t)\) is approximated via the posterior mean estimate $\nabla_{\x_t}\log p(\y^{\delta}|\hat{\x}_0)$ for $\hat{\x}_0 := \Ed[\x_0|\x_t]$. For example \cite{chung2023diffusion}, in the Gaussian-noise case with variance \(\sigma^2\), one obtains under an \(\ell^2\)-norm
\begin{equation}\label{eq:likelihood_grad}
\begin{aligned}
    \nabla_{\x_t} \log p(\y^{\delta}|\hat\x_0(\x_t)) &= \frac{\partial \hat{\x}_0\left(\x_t\right)}{\partial \x_t} \nabla_{\hat{\x}_0} \log p(\y^{\delta}|\hat\x_0) \\
     &= - \frac{1}{\sigma^2}  \frac{\partial \hat{\x}_0\left(\x_t\right)}{\partial \x_t} 
     \nabla_{\hat{\x}_0}\left\|\y^{\delta}-\mathcal{A}\left(\hat{\x}_0\right)\right\|_2^2 \\ 
     &= - \frac{2}{\sigma^2} \frac{\partial \hat{\x}_0\left(\x_t\right)}{\partial \x_t} \left(\partial \Ac \right)^*_{\hat{\x}_0}\left(\mathcal{A}\left(\hat{\x}_0\right) - \y^{\delta}\right), \end{aligned}  
\end{equation}  
where \((\partial\Ac)^*_{\hat\x_0}\) is the adjoint of the Fréchet derivative at \(\hat\x_0\).

\section{Main Method}
\label{sec:4} In this section, we will discuss a novel approach to solve USCT by leveraging \textit{adjoint neural operators}, seamlessly integrated with the  \textit{conditional consistency model} as data prior. Here, we will first introduce how to utilize the neural operator to approximate the adjoint-based gradient in \cref{sec:4.1}. 

\subsection{Adjoint Operator Learning via Multi-grid Neural Operator}\label{sec:4.1} 
In USCT, reconstructing the sound-speed distribution $\X_0(\mathbf{r})$ from the noisy measurement $\y^\delta$ constitutes a PDE-constrained optimization problem. The formulation of this inverse problem can be expressed explicitly as follows:
\begin{equation}
\begin{aligned}
\min_{\X_0}  \quad \mathcal{T}(\X_0) 
 &= \sum_{m}\sum_{n}\|\Y_n(\mathbf{r}_m) - \y_{m,n}^\delta \|_2^2, \\ 
\textit{s.t.} \quad \nabla^2 \Y_n(\mathbf{r}) +& \frac{\omega^2}{\X_0(\mathbf{r})^2}\Y_n(\mathbf{r}) = -\rho_n(\mathbf{r}).
\end{aligned}
\end{equation}

Gradient-based optimization methods, particularly the adjoint-state methods \cite{plessix2006review, 
taillandier2009first}, provide efficient mechanisms for computing gradients in PDE-constrained optimization problems. This method introduces an auxiliary adjoint field $\boldsymbol{\Lambda}_n(\mathbf{r})$, governed by the following adjoint equation:
\begin{equation}\label{eq:adjoint}
\nabla^2 \boldsymbol{\Lambda}_n(\mathbf{r}) + \frac{\omega^2}{\X_0(\mathbf{r})^2}\boldsymbol{\Lambda}_n(\mathbf{r}) = \sum_{m=1}^{M} \overline{\left(\Y_n(\mathbf{r}_m) - \y_{m,n}^\delta\right)} \rho_n(\mathbf{r}).
\end{equation}

A critical observation is that the adjoint equation shares the same structure as the forward equation. Specifically, the Helmholtz operator $\left[\nabla^2  + \frac{\omega^2}{\X_0(\mathbf{r})^2}\right](\cdot)$ is self-adjoint under suitable boundary conditions, such as the absorbing boundary conditions applied here to emulate wave radiation into an infinite domain. Consequently, solving the adjoint equation equates to solving the forward equation with adjusted source terms. By exploiting self‐adjointness and linearity, one shows that
\begin{equation}\label{eq:linear_adjoint}
\begin{aligned}
\boldsymbol{\Lambda}_n(\mathbf{r}) = -\sum_{m=1}^{M} \overline{\left(\Y_n(\mathbf{r}_m) - \y_{m,n}^\delta\right)} \Y_m(\mathbf{r}).
\end{aligned}
\end{equation}

Crucially, no separate discretization or solve of the adjoint equation~\cref{eq:adjoint} is required. All computations reduce to a batch of forward Helmholtz solves solutions $\{\Y_n\}^{N}_{n = 1}$ for all source points. This forward‐only formulation directly enables the efficient gradient updates for USCT without introducing new adjoint PDE solvers.

\begin{figure}[htp!]
    \centering
    \includegraphics[width=1.0\linewidth]{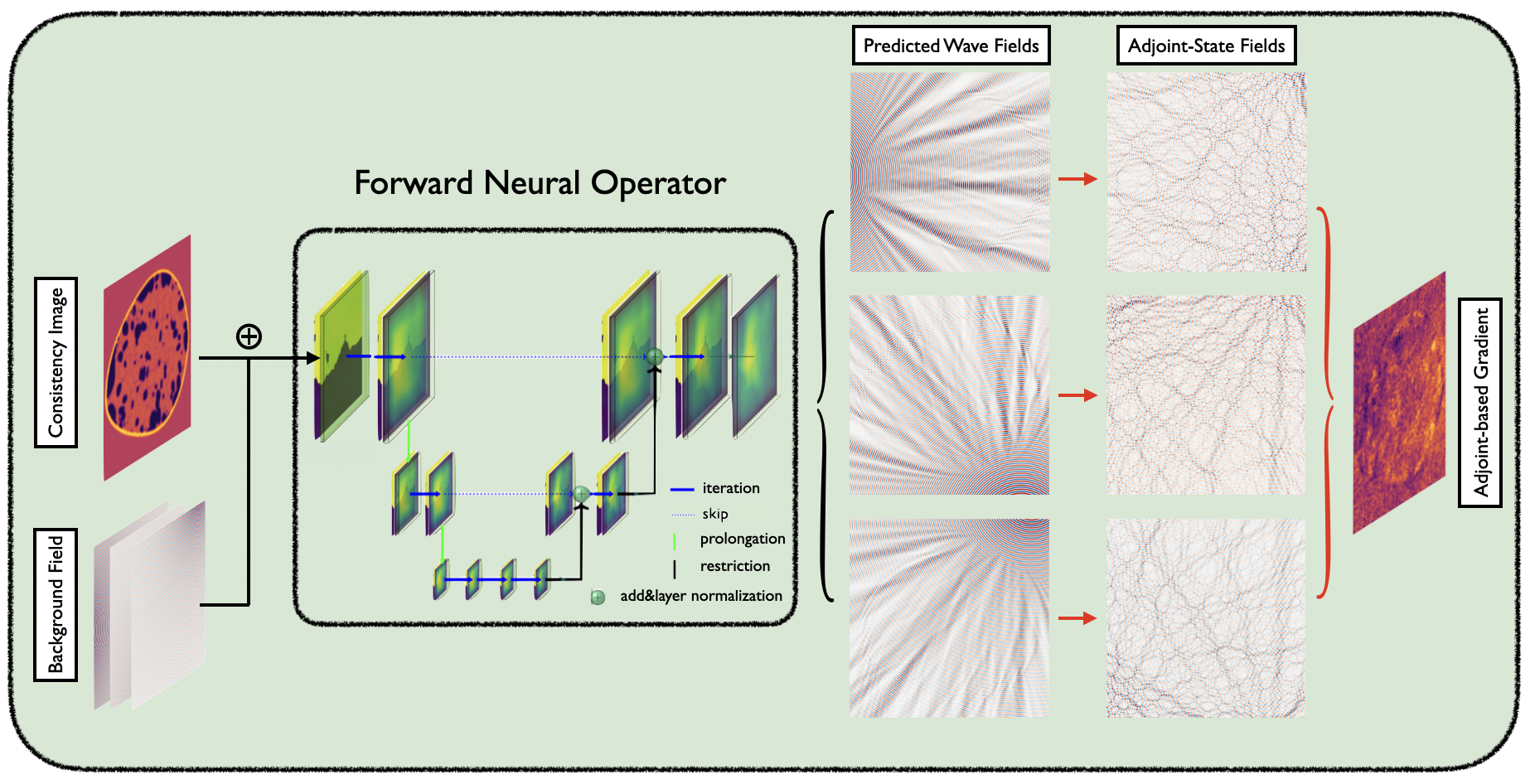}
    \caption{Overview of the \textit{Adjoint Neural Operator} framework. A consistency image generated from CM and a background wavefield pre-computed from the homogeneous Helmholtz equation are combined as input channels. The MgNO predicts the forward wavefields. Adjoint-state fields are computed to form adjoint-based gradients, enabling efficient PDE-constrained optimization.}
    \label{fig:Neural_Adjoint_Optimization}
\end{figure}

\textbf{Multi-grid Neural Operator.} Sampling requires hundreds of forward/adjoint Helmholtz solves with different $\X_0$ and source $\rho$, and numerical schemes such as the \emph{Convergent Born Series} \cite{osnabrugge_convergent_2016} are therefore the computational bottleneck. To accelerate the sampling, we use neural operator as a surrogate instead of numerical methods. In this work, we adapt a MgNO \cite{he2023mgno}, denoted $\widetilde{\mathcal{G}}_{\Theta_{\rm Mg}}$, specifically architected for solving the Helmholtz equation. The network takes a sound-speed map $\X_0$ and a source term $\rho$ to return the wavefield approximation: $\widetilde{\mathcal G}_{\Theta_{\mathrm{Mg}}}\bigl(\X_0 ;\rho_n \bigr)\;\approx\;\Y_n$.

We implement $\widetilde{\mathcal G}_{\Theta_{\mathrm{Mg}}}$ as a learnable iteration block following a multigrid V-cycle structure \cite{hackbusch2013multi}, where the restriction, prolongation, and smoothing operations are parameterized by convolutional neural blocks. This architecture allows the operator to perform multi-scale error correction through hierarchical downsampling and upsampling, as standard in multigrid methods, while adapting to spatial variations in the medium. Once $\widetilde{\mathcal{G}}_{\Theta^{*}_{\rm Mg}}$ is well-trained, as shown in \cref{fig:Neural_Adjoint_Optimization}, all forward Helmholtz solves within the gradient computation are replaced by network evaluations. Then, the gradient of the loss can be approximated by:
\begin{equation}\label{eqs:approx_gradient}
\frac{\partial \mathcal{T}}{\partial \X_0}(\mathbf{r}) 
\approx \frac{2\,\omega^2}{\X_0(\mathbf{r})^3}
\sum_{n=1}^N \sum_{m=1}^M 
\overline{\bigl(\widetilde{\Y}_n(\mathbf{r}_m) - \y_{m,n}^\delta\bigr)}\,
\widetilde{\Y}_m(\mathbf{r})\,
\widetilde{\Y}_n(\mathbf{r}),
\end{equation}
where \(\widetilde{\Y}_m := \widetilde{\mathcal{G}}_{\Theta^{*}_{\rm Mg}}\bigl(\X_0\,;\rho_m \bigr)\) and $\widetilde{\Y}_n  := \widetilde{\mathcal{G}}_{\Theta^{*}_{\rm Mg}}\bigl(\X_0\,;\rho_n \bigr)$ denote the surrogate solutions for source \(\rho_m\) and \(\rho_n\), respectively. By batching all evaluations of $\{ \widetilde{\Y}_n\}_{n = 1}^{N}$ on modern GPU hardware, the computational cost becomes dominated by neural network inference, resulting in orders‐of‐magnitude speedups compared to traditional PDE solvers. 

\textbf{Online Training.} To streamline the training of $\widetilde{\mathcal{G}}_{\Theta_{\rm Mg}}$, we introduce the online BCBS to enable GPU-accelerated parallel computation of CBS across multiple sources and data samples. CBS was originally proposed to solve the Helmholtz equation in arbitrary strong scattering media by constructing a convergent preconditioner, but its sample-by-sample iterative nature results in slow convergence rates in practice. Our BCBS effectively unrolls these iterations as batched tensor operations over the GPU, computing updates for multiple independent source terms and sound-speed environments simultaneously. This structural modification exploits the fact that the Green's operator remains constant across samples, thus enabling the rapid, memory-efficient computation of wavefield supervision targets $\{\Y_n\}^{N}_{n=1}$ on-the-fly during training. This directly allows mini-batch stochastic gradient descent to be applied over highly diverse, dynamically generated sound-speed fields. We summarize the pipeline in \cref{alg:batch_cbs_train_NO}. 

\begin{algorithm}[t]
\caption{Training Neural Operator using BCBS from CM}\label{alg:batch_cbs_train_NO}
\KwIn{
(1) \textit{CBS Parameters}: Batch‐based CBS solver $\operatorname{BCBS}(\overline{\overline{\X}}_0\, ; \overline{\overline{\rho}}_{\rm full}, \omega, D)$, full source terms $\overline{\overline{\rho}}_{\rm full}\in\mathcal{X}(\Omega)^{1\times N}$, angular frequency $\omega$, maximum iterations $D$; \\
(2) \textit{Network Parameters}: pretrained consistency model $f_{\Theta^{*}_{\rm CM}}(\x_t, t)$ for $t\in[\varepsilon,T]$, untrained neural operator $\widetilde{\mathcal{G}}_{\Theta_{\rm Mg}}(\X_0\,; \overline{\overline{\rho}}_{\rm full})$; \\
(3) \textit{Training Settings}: training batch size $N_0$, training epochs $E$, optimizer $\mathrm{Opt}(\cdot)$
}
\KwOut{Trained neural operator $\widetilde{\mathcal{G}}_{\Theta^{*}_{\rm Mg}}(\X_0 ; \overline{\overline{\rho}}_{\rm full})$}

\For{$e = 1$ \KwTo $E$}{
    // \texttt{Sample from the consistency model}\;
    $\overline{\overline{\x}}_0 \in \mathbb{C}^{N_0\times 1\times H \times W}$ $\gets$ Sample a batch of $N_0$ fields from $f_{\Theta^{*}_{\rm CM}}(\overline{\overline{\x}}_t, t)$\;
    $\overline{\overline{\X}}_0 \in \mathcal{X}(\Omega)^{N_0 \times 1}$ $\gets$ Interpolate $\x_0$ from image to physical domain for $\widetilde{\mathcal{G}}_{\Theta_{\rm Mg}}(\cdot \, ; \overline{\overline{\rho}}_{\rm full})$\; 
    
    // \texttt{Generate CBS supervision targets on-the-fly }\;
    $ \overline{\overline{\Y}} \in \mathcal{Y}(\Omega)^{N_0 \times N} \gets$ $\operatorname{BCBS}(\overline{\overline{\X}}_0\, ; \overline{\overline{\rho}}_{\rm full}, \omega, D)$\; 
    
    // \texttt{Neural operator training }\;
    $\mathcal{T}\gets$ Loss function $\mathcal{T} := \bigl\|\widetilde{\mathcal{G}}_{\Theta_{\rm Mg}}\!\left(\overline{\overline{\X}}_{0};\, \overline{\overline{\rho}}_{\rm full}\right)-\overline{\overline{\Y}}\bigr\|_{2}^{2}$\;
    $\Theta_{\rm Mg}\leftarrow \mathrm{Opt} \bigl(\Theta_{\rm Mg},\nabla_{\Theta_{\rm Mg}}\mathcal{T}\bigr)$\; 
}
\end{algorithm}

\subsection{Conditional Consistency Model with Adjoint Neural Operator}\label{sec:4.2}
The fundamental idea of a conditional CM is to introduce conditions explicitly to guide the self-consistent multi-step refinement of reconstructions. Following the approach outlined in \cite{zhao2024cosign}, we prioritize optimizing the following direct reconstruction loss:
\begin{equation}
\mathcal{T}_{\rm recon} := d\bigl(f_{\Theta_{\rm CC}}(\x_t, \y^{\delta}, t) , \x_0\bigr),
\end{equation}
rather than the consistency loss defined as \cref{eq:consistency_loss}. To balance the direct reconstruction loss with the consistency constraint, we can introduce an additional control block (also known as ControlNet \cite{zhang2023adding}) over the frozen CM backbone for guiding the multi-step conditional generation. The architecture and initialization of the control block follow \cite{zhao2024cosign, zhang2023adding}, where conditions are embedded into the consistency model through zero-convolution adapters. The structure of the control block and consistency model is shown in \cref{fig:framework}. In this way, the multi-step generation ability can be inherited from the pretrained frozen CM, while ensuring the reconstruction consistent with the measurements. Once the optimal $\Theta^{*}_{\rm CC}$ is trained, the multi-step conditional sampling is to repeatedly applies
\begin{equation}
\begin{aligned}
     \x_0 &\leftarrow f_{\Theta^{*}_{\rm CC}}\bigl(\x_{\tau_{n+1}},\y^{\delta},\tau_{n+1}\bigr), \\
     \x_{\tau_n}  &\leftarrow \text{Forward SDE}(\x_0,\tau_n), \quad n=S-1,\dots,0,
\end{aligned}
\end{equation}
for a given sampling-time sequence $\varepsilon = \tau_0 < \tau_1 < \cdots < \tau_{S} = T$, with the initialization $\x_{\tau_{S}}$ from the Gaussian prior distribution. 

Note that ControlNet is originally designed to handle image-domain conditions, requiring the measurements $\y^{\delta}$ mapped into the image-domain as the input for the control block. \cite{zhao2024cosign} employ the pseudo-inverse operator $\mathcal{A}^\dagger$ to generate an initial reconstruction for linear inverse problems, whereas they feed the resized measurements directly as the condition for nonlinear scenarios. For ill-posed USCT, as detailed in \cref{sec:1.1}, traditional iterative approaches often become trapped in local minima and fail to produce satisfactory pre-reconstruction for conditioning the control block. To ensure both speed and fidelity, we introduce a supervised-based network to directly approximate the inverse mapping $\y^{\delta}\mapsto\x_0$. Thereafter, the conditional CM further refines the pre-reconstruction that better aligns with the prior data distribution \footnote{Empirically, we observe that the pre-reconstruction fidelity correlates positively with the final performance of the conditional CM framework, as demonstrated in our ablation study (See \cref{sec:7.1}).}. 

\textbf{Adjoint Neural Operator for Measurement-Informed Guidance.} While the conditional CM provides an initial data-driven mapping, achieving high-fidelity reconstructions requires strict alignment with the observed acoustic measurements. Rather than characterizing any valid positive sound speed as "physics-informed," this step is strictly a \textit{measurement-informed} or \textit{data-informed} projection. During reverse sampling, optimization algorithms iteratively update the prior sample $\x_0$ to push it onto the measurement-consistent manifold $\mathcal{M} := \{ \z \, \big| \,  \| \mathcal{A}(\boldsymbol{z}) - \boldsymbol{y}^{\delta} \|_2^2 \leq \epsilon^2 \}$. Because we are explicitly minimizing this least-squares residual (analogous to solving \cref{eq:inverse}), the final point estimates correspond to the measurement-consistent projection (a refined conditional expectation of the posterior), rather than a purely unconditional generation.

\begin{algorithm}[t]
\caption{Diff-ANO: Conditional sampling with adjoint neural operator for USCT}
\label{alg:sampling}
\KwIn{
(1) \textit{Neural Operator Parameters}: observed wavefield data $\y^\delta \in \mathbb{R}^{M \times N}$, forward neural operators $\widetilde{\mathcal{G}}_{\Theta^*_{\rm Mg}}(\X_0 ; \rho_n)$ for $n = 1, \dots, N$, receiver points $\{\mathbf{r}_m\}_{m=1}^{M}$; \\
(2) \textit{Consistency Model Parameters}: conditional consistency model $f_{\Theta^{*}_{\rm CC}}(\x_t, \y, t)$ for $t \in [\varepsilon , T]$, sequence of sampling-time points $ \varepsilon = \tau_0 < \dots < \tau_S = T$
}
\KwOut{$\x_0$}

\textbf{Initialize:} $\x_0 \gets f_{\Theta^{*}_{\rm CC}}(\x_T, \y^\delta, T)$ where $\x_T$ follows the prior distribution\;

\For{$s \gets S-1$ \KwTo $1$}{
    // \texttt{Perform neural adjoint optimization for guidance}\;
    $\X_0 \gets$ Interpolate $\x_0$ from image-domain to physics-domain for $\widetilde{\mathcal{G}}_{\Theta^*_{\rm Mg}}(\cdot ; \rho_n)$\;
    
    \If{\textit{neural adjoint optimization}}{
        $\{\widetilde{\Y}_n\}^{N}_{n = 1} \gets \widetilde{\mathcal{G}}_{\Theta^*_{\rm Mg}}(\X_0 ; \rho_n)$ for $n = 1, \dots, N$\;
        $\widetilde{\boldsymbol{\Lambda}}_n \gets -\sum_{m=1}^{M}\overline{ \left(\widetilde{\Y}_n(\mathbf{r}_m) - \y^\delta_{m,n}\right)} \widetilde{\Y}_m$\;
        Update $\X_0$ from the adjoint-based gradient $\nabla_{\X_0}\mathcal{T}(\X_0) = - \frac{2\omega^2}{\X_0^3}\sum_{n=1}^{N} \widetilde{\boldsymbol{\Lambda}}_n \widetilde{\Y}_n$\;
    }
    
    $\x_0 \gets$ Interpolate $\X_0$ from physics-domain to image-domain for $f_{\Theta^{*}_{\rm CC}}(\cdot , \y^\delta, \tau_s)$\;
    
    // \texttt{Multi-step consistency mapping}\;
    $\x_{\tau_s} \gets$ Sample from the forward SDE at time $\tau_s$ with the initial $\x_0$\;
    $\x_0 \gets f_{\Theta^{*}_{\rm CC}}(\x_{\tau_s}, \y^\delta, \tau_s)$\;
}
\Return $\x_0$\;
\end{algorithm}

Furthermore, in Algorithm \ref{alg:sampling}, the adjoint-based gradient step can be applied either once or iteratively multiple times until convergence at each consistency step, depending on the desired strictness of the measurement constraint. For nonlinear inverse problems, where closed-form solutions are unavailable, gradient descent or its variants are typically employed to solve it. In PDE-governed scenarios, however, adjoint methods are preferred to compute the data-fidelity gradient $\nabla_{\boldsymbol{z}}\| \mathcal{A}(\boldsymbol{z}) - \boldsymbol{y}^{\delta} \|_2^2$ rather than formulate the explicit Jacobian $(\partial \mathcal{A})_{\boldsymbol{z}}$. This strategy relies well-established numerical solvers for both forward and adjoint PDEs, but introduces two principal problems:
\begin{enumerate}
    \item \textit{Continuous-Domain to Discrete-Domain.} 
    When solving PDE-based inverse problems within the consistency model framework, discretization of the PDE inevitably introduces numerical approximation errors during the multi-step sampling, resulting in the reconstruction fidelity degradation. It is necessary to bridge the gap between the continuous physics-domain where the governed PDEs are naturally formulated, and the discrete image-domain where the consistency model acts.  

    \item \textit{PDE-Solver Bottleneck.} 
    Each evaluation of the data-fidelity gradient requires solving both the forward and adjoint PDEs under multiple boundary conditions. Empirically, these numerical PDE solves dominate the computational budget and become the bottleneck of the conditional sampling process, whereas the consistency model can leverage efficient GPU acceleration. Achieving real-time performance thus demands the acceleration of PDE-based gradient estimates. 
\end{enumerate}
In USCT, to handle the aforementioned issues, we propose to incorporate the \textit{adjoint neural operator}, as detailed in \cref{sec:4.1},  to impose physics-informed constraints into the multi-step conditional sampling process. By utilizing the pretrained neural operators, the computational cost in \cref{eqs:approx_gradient} becomes dominated by network inference, resulting in orders-of-magnitude speedups compared to traditional CBS solvers. Besides, due to the inherent discretization-invariance \cite{boulle2024mathematical} of the neural operator, the common interpolation strategy can be adopted to transform the sample between physics-domain and image-domain. Here, the overall sampling algorithm is shown in \cref{alg:sampling}.

\begin{figure}[htp!]
    \centering
        \begin{minipage}{0.475\textwidth}
        \centering
        \includegraphics[width=\textwidth]{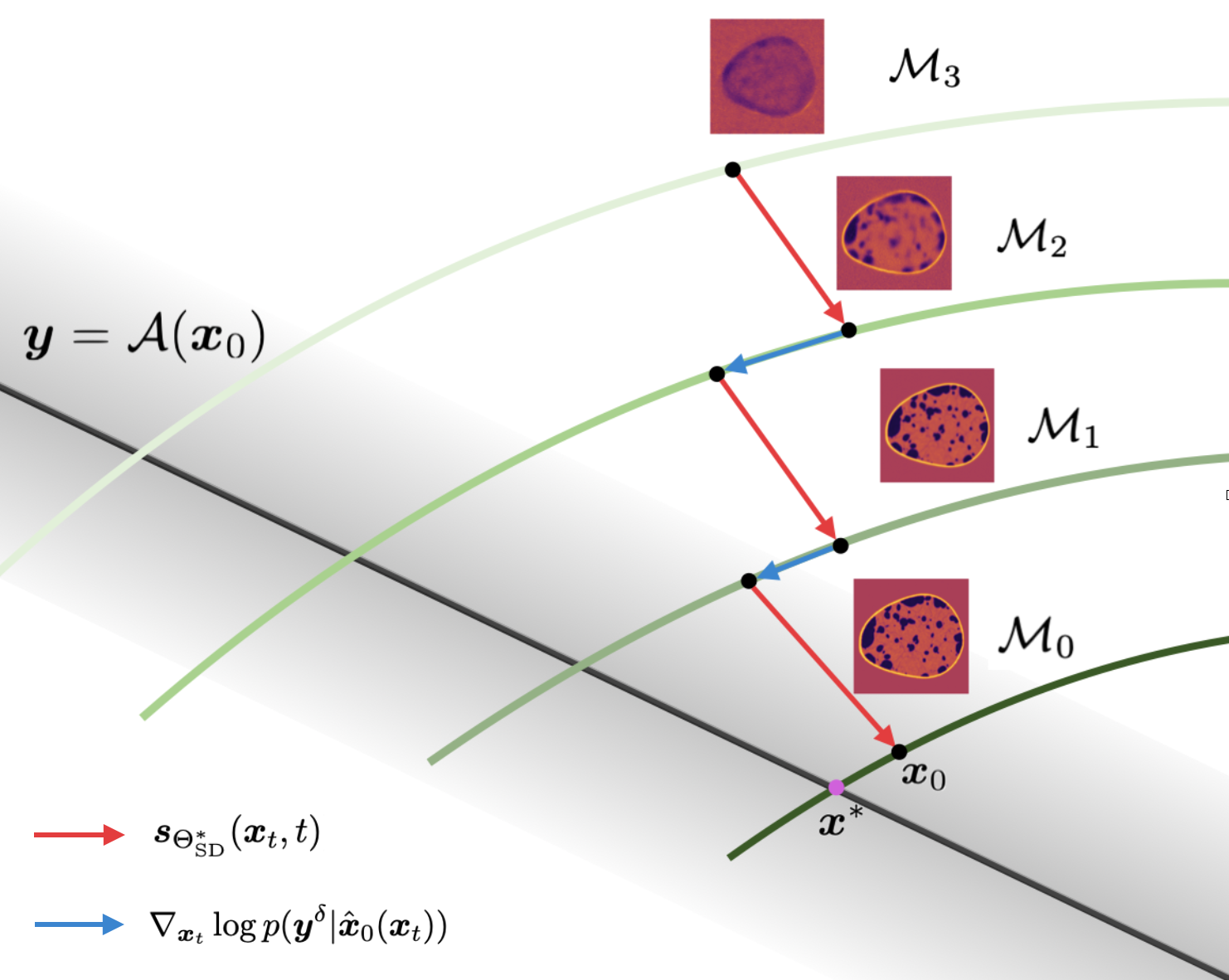}
    \end{minipage}%
    \hfill
    \begin{minipage}{0.47\textwidth}
        \centering
        \includegraphics[width=\textwidth]{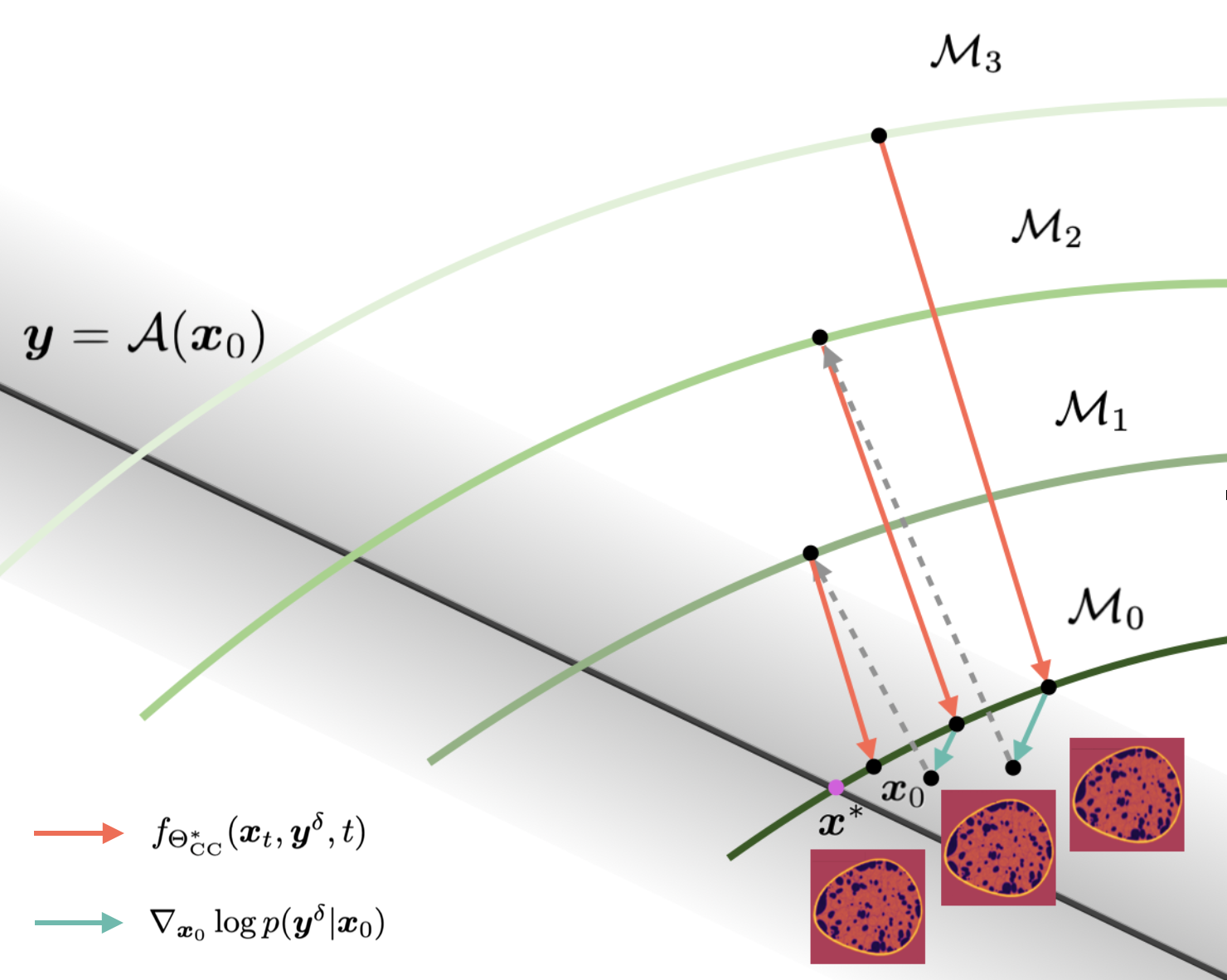}
    \end{minipage}
    
    \caption{Comparison between the trajectory of DPS \cite{chung2022diffusion} and Ours. The gray region denotes the measurement-consistent manifold  $\mathcal{M} := \{ \x_0 \, \big| \,  \| \mathcal{A}(\x_0) - \boldsymbol{y}^{\delta} \|_2^2 \leq \epsilon^2 \}$, and the green curves with varying saturation represent the distribution $\mathcal{M}_t$ of noised samples $\boldsymbol{x}_t$. The ground truth is denoted as $\boldsymbol{x}^{*}$. \textbf{Left:} DPS combines the score-based reverse sampling (red arrows) with the gradient guidance (blue arrows) updated on the Tweedie approximation $\hat{\x}_0 := \mathbb{E}[\boldsymbol{x}_0 | \boldsymbol{x}_t]$; \textbf{Right:} Ours incorporates the conditional consistency sampling (orange arrows) and the gradient guidance (cyan arrows)  updated on the clean sample $\boldsymbol{x}_0$.}
    \label{fig:trajectory_illustration}
\end{figure}

\textbf{Remarks.} It should be noted that the adjoint neural-operator based on $\widetilde{\mathcal{G}}_{\Theta^{*}_{\rm Mg}}$ cannot be directly applied to methods such as DPS \cite{chung2022diffusion} or CBS-based gradient descent \cite{osnabrugge_convergent_2016}. Here, we illustrate this by comparing the sampling trajectories of DPS and conditional CM in \cref{fig:trajectory_illustration}. In our framework, the application of the trained neural operator $\widetilde{\mathcal{G}}_{\Theta^{*}_{\rm Mg}}$ critically depends on its generalization capability, where $\widetilde{\mathcal{G}}_{\Theta^{*}_{\rm Mg}}$ is merely trained on the clean data manifold \(\mathcal{M}_0\). On the one hand, since the conditional CM maps measurements onto \(\mathcal{M}_0\), the operator $\widetilde{\mathcal{G}}_{\Theta^{*}_{\rm Mg}}$ can be directly plugged in to estimate the gradient guidance. On the other hand, DPS requires computing the gradient over the augmented manifold
\[
\overline{\mathcal{M}}_0 := \bigl\{\hat{\x}_0 \mid \hat{\x}_0 = \mathbb{E}[\x_0| \x_t],\, \x_t \in \mathcal{M}_t,\ \forall t \in [0,T]\bigr\}.
\]
Since \(\mathcal{M}_0 \subset \overline{\mathcal{M}}_0\), directly using $\widetilde{\mathcal{G}}_{\Theta^{*}_{\rm Mg}}$ trained on \(\mathcal{M}_0\) degrades its generalization ability in adjoint-based optimization, leading to inaccurate gradient estimates during DPS. Furthermore, to handle the ill-posed nature of USCT, the conditional CM provides an appropriate initial reconstruction for optimization, whereas DPS lacks such an initial prior.

\section{Implementations}\label{sec:5}

\subsection{Dataset Collection and Measurement Configuration}\label{sec:5.1}
Our numerical experiments utilize the phantom dataset provided by the \textit{OpenBreastUS} dataset \cite{zeng2025openbreastusbenchmarkingneuraloperators}, a comprehensive, anatomically realistic USCT resource for benchmarking neural wave equation solvers. \textit{OpenBreastUS} consists of 8000 breast phantoms, categorized into four distinct groups based on breast density characteristics: heterogeneous (HET), fibroglandular (FIB), all fatty (FAT), and extremely dense (EXD). For our experiments, we exclusively select the FIB and EXD subsets, comprising 2,700/300 and 1,800/200 training-testing samples respectively. All simulations assume a uniform background sound-speed value of 1500 m/s, with regions of interest (ROI) exhibiting heterogeneous sound-speed distributions ranging between 1408 m/s and 1595 m/s. To further validate cross-anatomical generalization, we also evaluate our method on the \textit{OpenPros} dataset \cite{wang2025openpros}, a limited-angle prostate USCT benchmark. While OpenBreast results are presented in detail throughout \cref{sec:6.1}--\cref{sec:6.5}, OpenPros results are summarized concisely in \cref{sec:6.6} with key quantitative metrics in \cref{tab:openpros_results}.

\begin{figure}[htp!]
    \centering
        \begin{minipage}{0.4\textwidth}
        \centering
        \includegraphics[width=\textwidth]{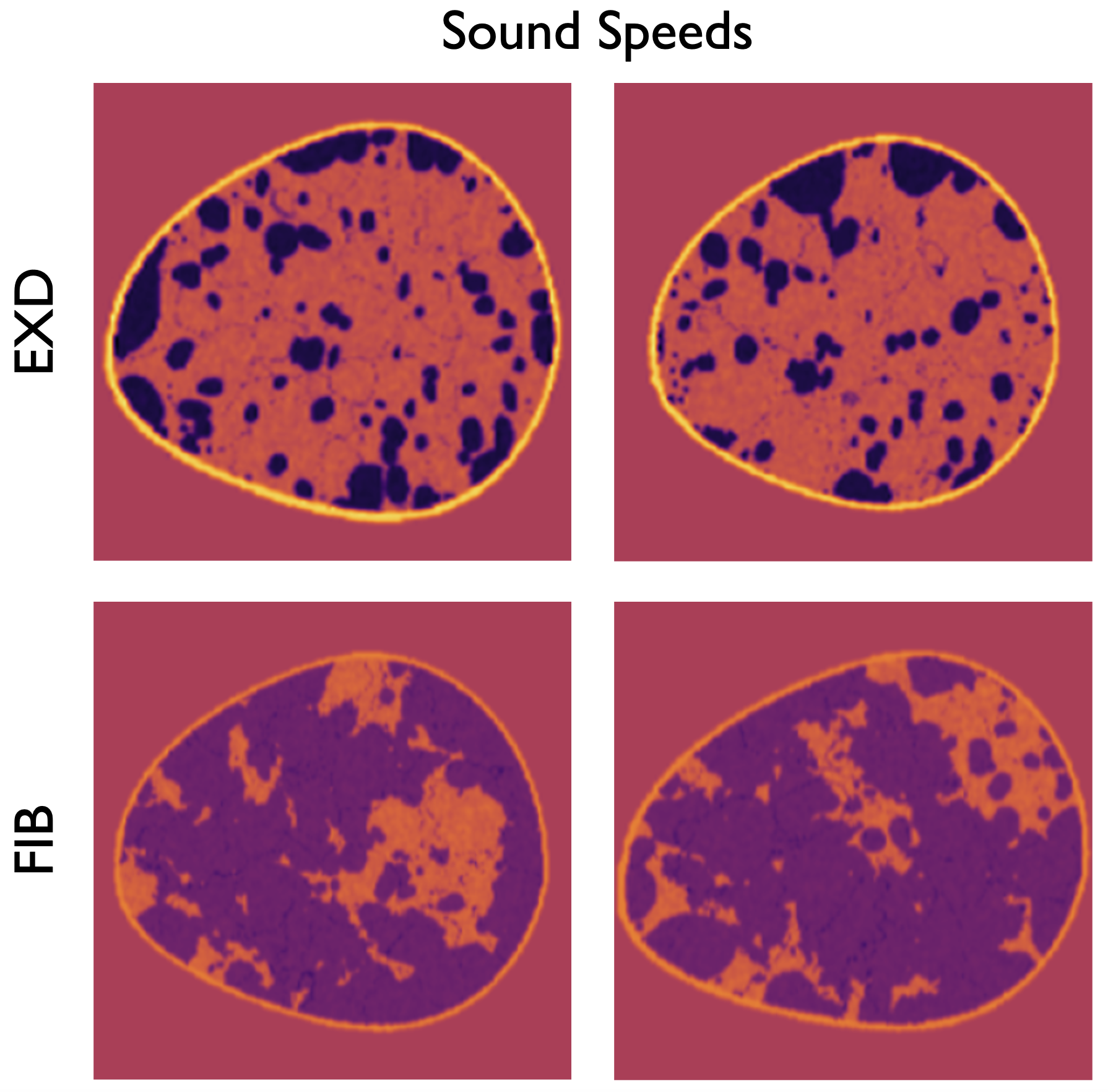}
    \end{minipage}
    \begin{minipage}{0.4\textwidth}
        \centering
        \includegraphics[width=\textwidth]{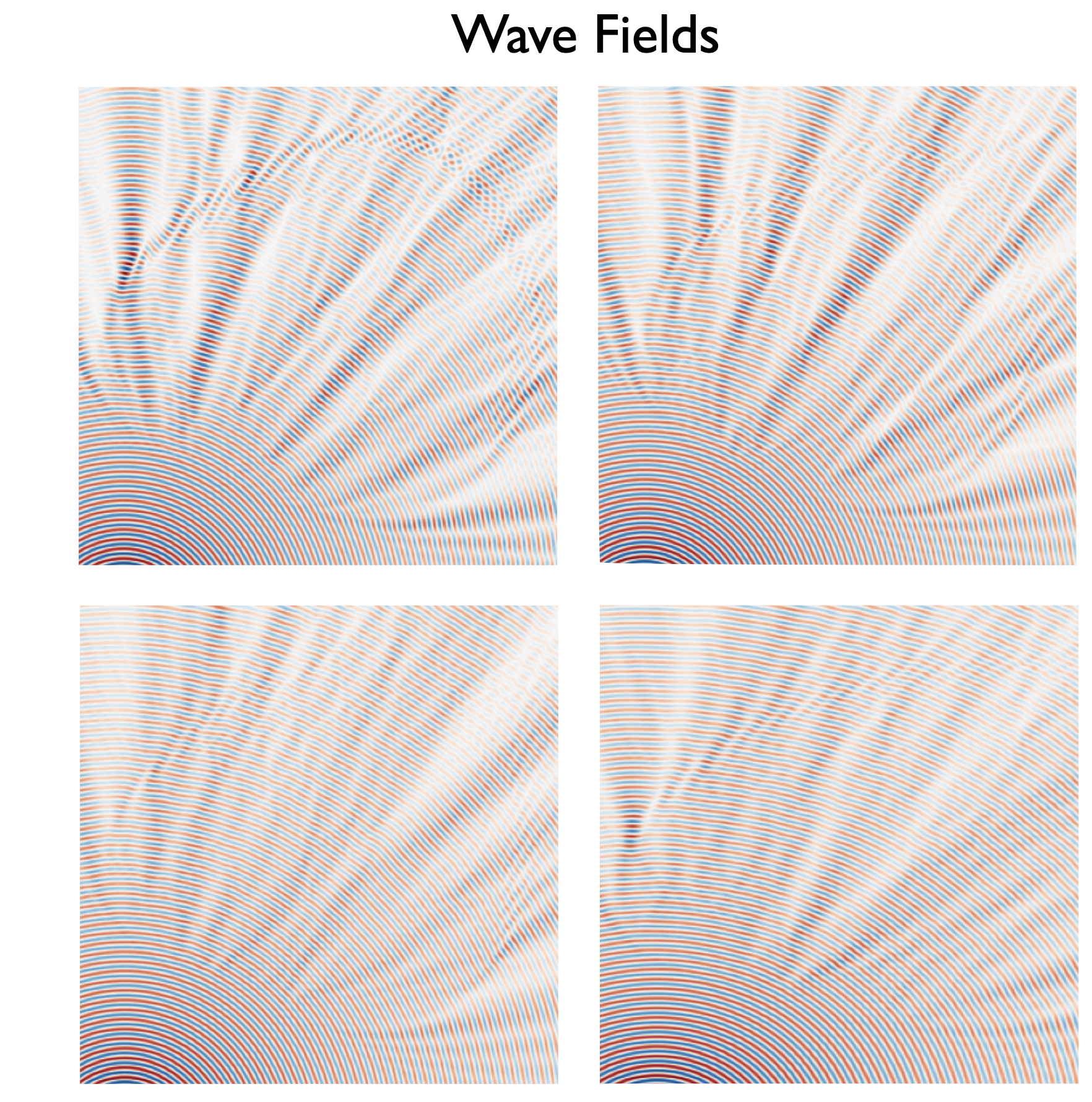}
    \end{minipage}
    \caption{The left panel shows heterogeneous sound–speed phantoms for two specific breast density types: \textbf{EXD} (top row) and  \textbf{FIB}(bottom row). The right panel displays the corresponding scattering wavefields in ROI obtained by solving the Helmholtz equation for a point source.}
    \label{fig:sample_illustration}
\end{figure}

\paragraph{Wavefield Generation.} Our experiments operate at a fixed frequency of 500 kHz, contrasting the multi-frequency approach employed in prior studies \cite{zeng2023neural}. In each measurement configuration, the observed data is collected from the resulting wavefields, which are simulated by leveraging CBS. Additionally, three different SNR levels (noise-free, 10dB, 5dB) are introduced into the observed data to evaluate the robustness of methods. Note that the speed samples and wavefields required for training the neural operator are generated on-the-fly by leveraging CM in \cref{sec:3.1} and BCBS in \cref{sec:4.1}, rather than utilizing precomputed and stored large datasets as \cite{zeng2025openbreastusbenchmarkingneuraloperators, wang2025openpros}. Here, four representative sound-speed samples and the corresponding wavefields are presented in \cref{fig:sample_illustration}. Crucially, FID-type samples feature weakly scattering medium, whereas EXD-type samples incorporate strongly scattering heterogeneities. 


\begin{figure}[htp!]
    \centering
    \includegraphics[width=\linewidth]{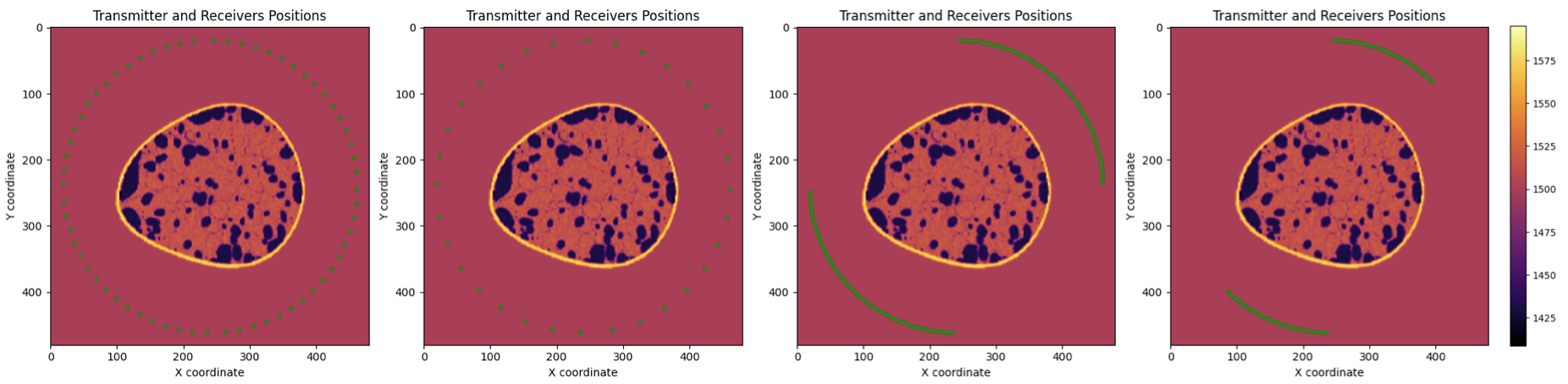}
    \caption{Measurement geometries considered in this work. Each panel shows a representative sound-speed distribution overlaid with the combined positions of transmitters and receivers (green dots).
    From left to right: (a) \textit{sparse-view} I; (b) \textit{sparse-view} II; (c) \textit{partial-view} I; (d) \textit{partial-view} II. 
    }
    \label{fig:measurement_settings}
\end{figure}

\paragraph{Measurement Configuration.} We adopt measurement configurations inspired by the physical settings presented in \cite{zeng2025openwaves}. Specifically, wavefields are simulated using parameters characteristic of a real annular USCT system, featuring $256$ transducers uniformly distributed around a $220$ mm diameter ring.  To rigorously evaluate the robustness of our approach under challenging scenarios, we introduce four under-sampled measurement configurations to systematically induce incomplete data conditions. These configurations, visualized in \cref{fig:measurement_settings}, are categorized into \textit{sparse-view} and \textit{partial-view} scenarios: (1) For sparse-view, we simulate $64$ source-receiver pairs uniformly distributed around the full ring, and $32$ source-receiver pairs for increasing the ill-posedness. (2) For partial-view, $64$ source-receiver pairs are uniformly distributed along a quarter-circle segment facing the ROI, and $32$ source-receiver pairs along an eighth-circle segment as well. These settings illustrate varying degrees of ill-posedness in terms of data incompleteness and angular coverage, designed to emulate real-world limitations commonly encountered in USCT, as detailed in \cref{sec:1.1}. 

\subsection{Network Architecture and Training}\label{sec:5.3}
The sound-speed distribution $\X_0$ is discretized on a $480 \times 480$ physics-domain grid with a spacing of $0.5$ mm. We adopt the mean-variance normalization to standardize it with $\mu = 1488.39$ and $\sigma = 27.53$, which are precomputed statistics from the training dataset.

\paragraph{MgNO Architecture} The MgNO architecture employs the following key configurations: Physical inputs $\Y^0_n$ (background wavefields) and $\X_0 $ (sound-speed distribution) are projected into latent states via $1 \times 1$ convolutional layers, sharing a unified feature dimension of $24$ channels. Each V-cycle block ${\mathcal{G}}_{\theta}$ executes iterative updates across seven resolution levels ($480, 239, 119, 59, 29, 14, 6$). Because these resolutions are non-dyadic, the coarsening strategy relies on adaptive average pooling (for downsampling) and bilinear interpolation (for upsampling) to seamlessly bridge fractional scale factors without losing spatial alignment. For adaptive convolution kernels, all $\mathcal{K}_h$ and $\mathcal{S}_h$ operators use $3\times 3$ convolutional filters in AdaConv layers, with MLP projections containing two hidden layers. We employ $6$ recurrent applications of the V-cycle iteration block.

\paragraph{Training Settings} In our implementation, we iteratively draw a batch of $N_0 = 32$ training samples $\{\x^{(i)}_0\}_{i=1}^{N_0}$ from the pretrained consistency model $f_{\Theta^{*}_{\rm CM}}$. Each sample $\x^{(i)}_0$, initially in the image-domain with a spatial resolution of $256 \times 256$, is transformed to $\X^{(i)}_0$ with a size of $480 \times 480$ through the bilinear interpolation. Then, we randomly select $N = 8$ source locations per batch, and the corresponding wavefields $\bigl\{\Y^{(i)}_1,\,\dots,\,\Y^{(i)}_N\bigr\}_{i=1}^{N_0}$ are produced from the batch-based CBS-solver to serve as the supervision targets. The MgNO is optimized by minimizing the empirical $\ell^2$‐loss, using the AdamW optimizer with an initial learning rate of $5 \times 10^{-4}$ and weight decay of $10^{-5}$. The training employs a OneCycleLR learning rate scheduler over $50$ epochs, following a cosine annealing strategy with $30\%$ warm-up phase. 

\section{Numerical Results}\label{sec:6}

\subsection{Reconstruction Results}\label{sec:6.1}
In this section, we compare the performance of the following  algorithms for USCT reconstruction using noise-5dB measurements in two distinct settings: sparse-$\mathrm{I}$ and partial-$\mathrm{I}$. For each scenario, we present the reconstructed results for two different sample types: EXD and FIB, as illustrated in \cref{sec:5.1}. The performance metrics, Peak Signal-to-Noise Ratio~(PSNR) and Structural Similarity Index Measure (SSIM), are provided for each reconstruction. 

\begin{figure}[htp!]
    \centering
    \includegraphics[width=1.0\linewidth]{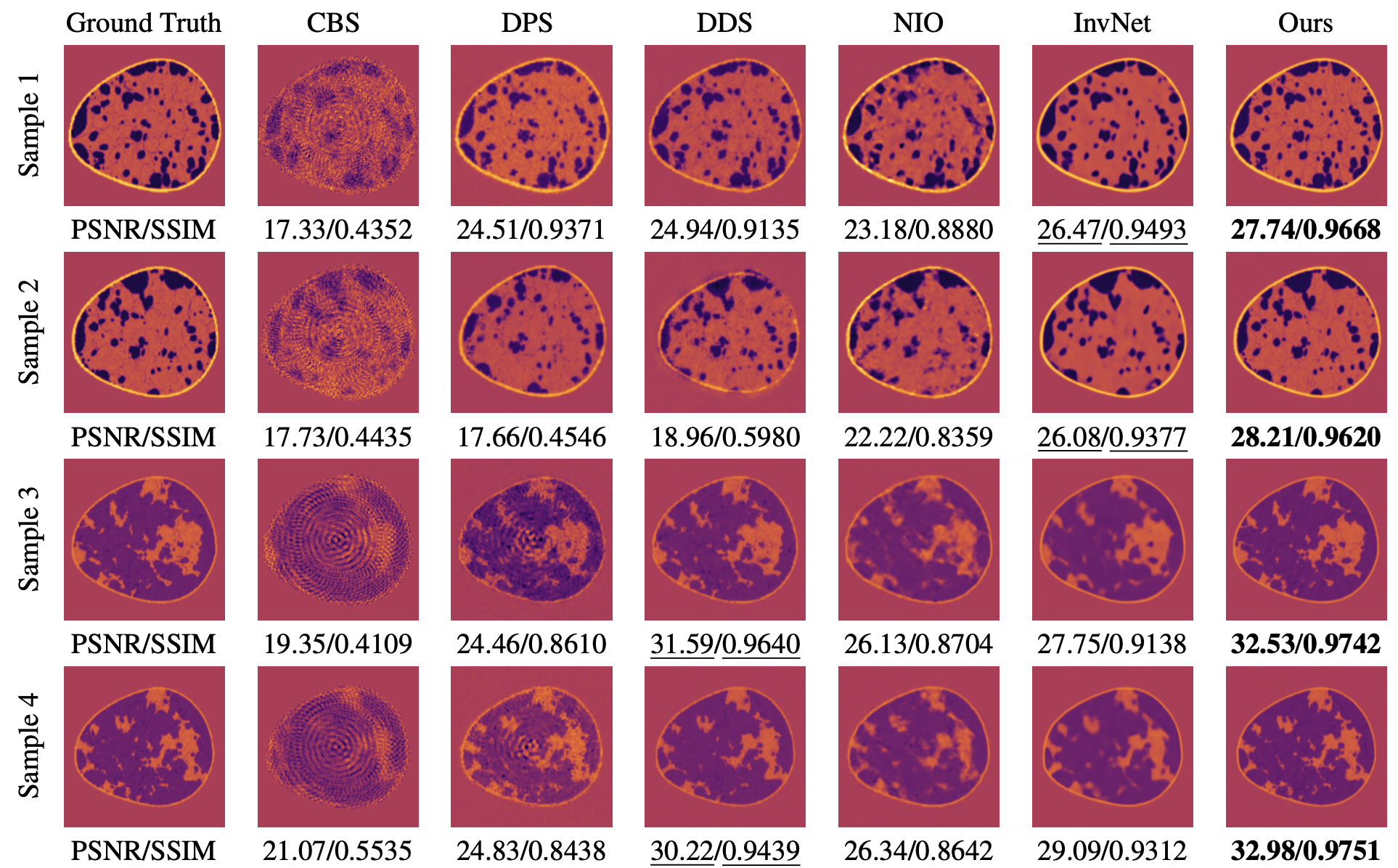}
    \caption{Reconstruction results under the sparse-view scenario with 5dB noise. The images represent the reconstructed results for different sample types, EXD (Samples 1 and 2) and FIB (Samples 3 and 4), and the PSNR and SSIM values are shown below.}
  \label{fig:sparse_recon}
\end{figure}

\begin{itemize}
    \item \textbf{CBS-Solver} \cite{osnabrugge2016convergent}: The CBS-Solver employs the Convergent Born Series method, to solve the forward Helmholtz equation with 64 source points. The optimization is performed using the adjoint-based method in \cref{sec:4.1} with the L-BFGS algorithm \cite{zhu1997algorithm}. In our experiments, the CBS-Solver uses $500$ inner iterations to ensure convergence for Helmholtz solutions and $30$ outer iterations for the L-BFGS optimization.

     
    \item \textbf{DPS} \cite{chung2022diffusion}:  The DPS method utilizes $1000$ discretized sampling steps over the sampling-time interval $[\varepsilon, T]$, and applies the batch-based CBS for faster optimization. At each sampling-time step, the BCBS is performed on Tweedie-denoised samples $\hat{\x}_0$ to solve the adjoint-based gradient and back-propagate it to the latent variable $\x_t$. Since the method requires the manual selection of the step size, we experimented  with multiple step sizes $\eta = 0.2, 0.1, 0.05$ and selected the best reconstruction result from these choices for comparison and analysis.

    \item \textbf{DDS} \cite{chung2023decomposed}: The DDS method extends DPS by decoupling the diffusion reverse sampling and gradient update steps. Conjugate Gradient descent \cite{shewchuk1994introduction} is applied on Tweedie-denoised samples, which eliminates the need for manual step-size selection. Additionally, DDIM sampling acceleration \cite{song2020denoising} is applied to expedite the posterior sampling process of DDS. In our experiments, we used a linear sampling strategy as DDIM with $50$ equidistant steps over the sampling-time interval $[\varepsilon, T]$. 

    \item \textbf{NIO} \cite{molinaro2023neural}: The Neural Inverse Operator (NIO) combines DeepONets and FNOs to approximate mappings from operators to functions. In our experiments, the convolutional layers in the branch network were adapted to match the spatial resolution requirements of USCT. The branch net employs a CNN comprising $9$ Conv2d layers to extract $512$ feature coefficients, which are subsequently projected onto $25$ basis functions through a linear transformation. The trunk net is implemented as an $8$-layer MLP with $100$ neurons per hidden layer. The FNO component utilizes $4$ Fourier layers, configured with $25$ Fourier modes and a channel width of $32$. 

    \item \textbf{Inversion-Net} \cite{zeng2021inversionnet3d}: Inversion-Net is a convolutional neural network featuring an encoder-decoder architecture designed for reconstructing 2D velocity distributions from seismic data. In our experiments, the network was trained on three distinct noise-level datasets: noise-free, 10dB and 5dB SNR conditions. The trained model was subsequently employed to predict 2D sound-speed maps from USCT observation data.

    \item \textbf{Ours}: Our approach leverages a pretrained CM as the sampling backbone and a pretrained MgNO as Helmholtz surrogates in the adjoint-based optimization. For all measurement scenarios, we use Inversion-Net  \cite{zeng2021inversionnet3d} as an inversion block to train the conditional consistency model as \cite{zhao2024cosign}. In our experiments, the sampling-time steps are chosen to be $\tau_1,\tau_2, \tau_3,\tau_4,\tau_5 = 0.1, 0.12, 0.14, 0.16, 0.18$  if the overall sampling-time interval is $[\varepsilon, T] = [0.001,1]$. 
\end{itemize}
In \cref{fig:sparse_recon}, the inherent ill-posedness of sparse-view reconstruction stems from insufficient measurement density. In CBS, traditional adjoint-based optimization without any regularization produces periodic oscillatory artifacts in the reconstructed images, result in the lowest PSNR and SSIM values among all methods. This underscores the necessity for appropriate regularization to suppress such artifacts in reconstructions. The DPS approach relies on a manually tuned step size for gradient update. It exhibits notable instability: excessively large step size leads to oscillatory artifacts as Sample 3, while excessively small step-size leads to insufficient physics constraints, yielding results consistent with the diffusion prior but significantly divergent from the ground truth as Sample 2. The DDS method effectively addresses the shortcomings of DPS by eliminating explicit step-size selection through CG and incorporating DDIM sampling acceleration. This provides more stable and significantly faster reconstruction compared to DPS. The effectiveness of diffusion-based priors in addressing the ill-posedness is clearly demonstrated in DPS and DDS. However, the supervised models exhibit distinct reconstruction artifacts: NIO introduces instability within homogeneous interior regions, resulting in noisy reconstructions. In contrast, Inversion-Net yields overly smooth outputs, and fails to recover subtle heterogeneities and textures. These supervised methods reveal the fundamental limitation of purely data-driven direct mapping approaches, constrained by sparse measurements without physics-informed optimization. Our proposed method integrates the benefits of direct mapping and iterative optimization, achieving faster, more stable, and superior reconstruction quality. 

In \cref{fig:partial_recon}, the limited angular coverage in this scenario introduces directional artifacts distinct from sparse-view oscillations, as evidenced in CBS results. This creates fundamentally different challenges that require implicit data completion through prior knowledge. For DPS, it again demonstrates sensitivity to the gradient update step size, particularly in reconstructing highly scattering EXD-type medium as Samples 1 and 2, though it shows relative stability for FIB-type reconstruction. DDS shows improved robustness over DPS in this scenario, but struggles to reconstruct EXD-type samples as well. Supervised learning approaches (NIO and Inversion-Net) exhibit similar limitations observed in the sparse-view scenario. Our proposed method, benefiting from consistency mapping and physics-informed optimization, consistently achieves the best PSNR/SSIM metrics and more accurately recovers image details, including subtle heterogeneities and internal textures.
\begin{figure}[htp!]
    \centering
    \includegraphics[width=1.0\linewidth]{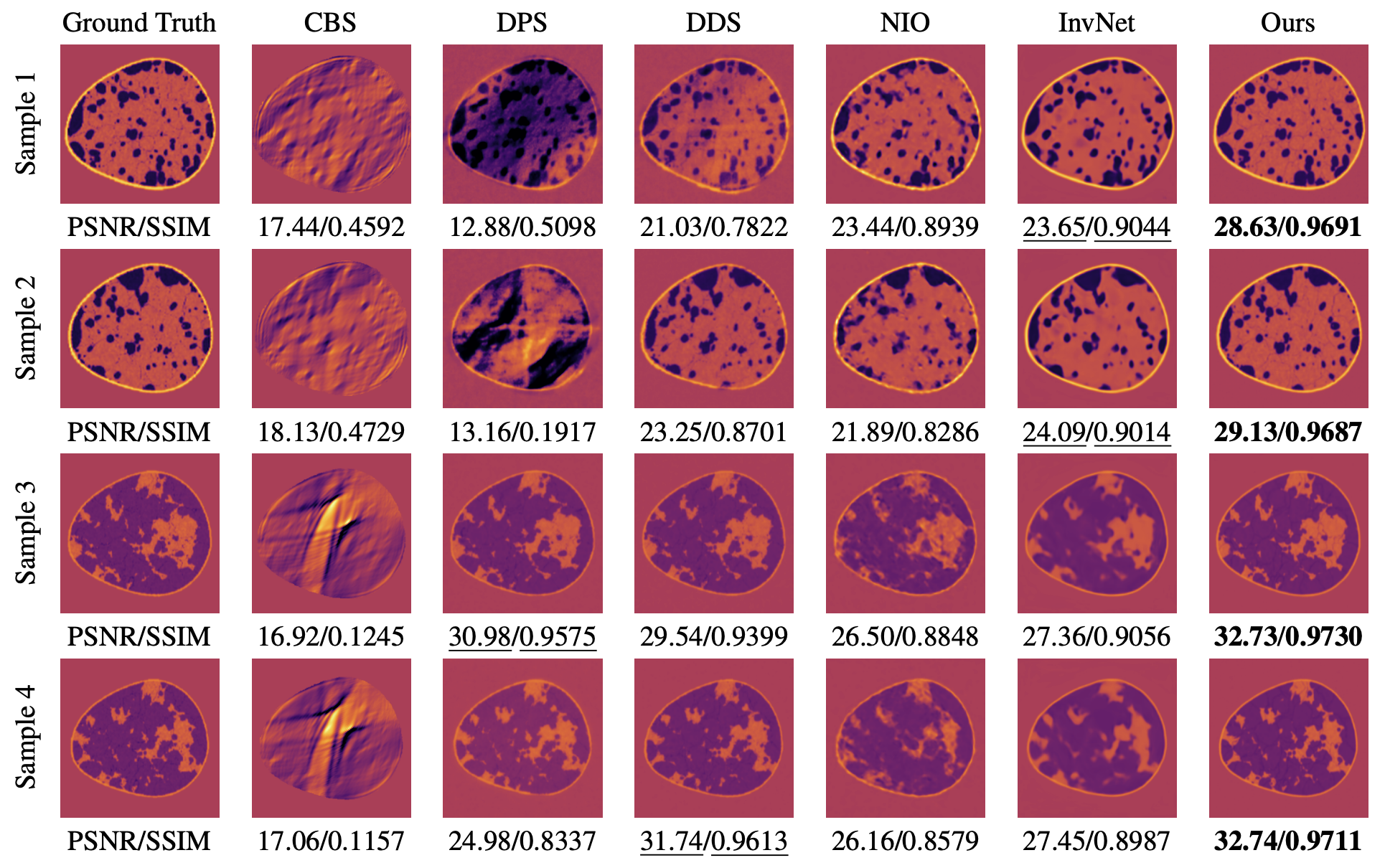}
    \caption{Reconstruction results under the partial-view scenario with 5dB noise. The images represent the reconstructed results for different sample types, EXD (Samples 1 and 2) and FIB (Samples 3 and 4), and the PSNR and SSIM values are shown below.}
  \label{fig:partial_recon}
\end{figure}

\textbf{Remarks.} From \cref{fig:sparse_recon} and \cref{fig:partial_recon}, the reconstruction performance for FIB-type samples consistently outperforms that for EXD-type samples in both scenarios. This can be attributed to weaker nonlinear scattering effects in FIB-type samples, as visualized in the right panel of \cref{fig:sample_illustration}, making their reconstruction closer to the linear inverse problem.

\subsection{Effect of the Optimization Step}\label{sec:6.2}
This subsection evaluates how the number of neural adjoint-based optimization steps influences reconstruction quality, particularly under challenging measurement scenarios with severe ill-posedness: sparse-$\mathrm{II}$ and partial-$\mathrm{II}$. We use the pretrained CM and the MgNO-based Helmholtz surrogates as detailed in \cref{sec:6.1}, progressively applying adjoint-based gradient guidance from later to earlier sampling-time steps ($\tau_5$ to $\tau_1$). 
\begin{figure}
    \centering
    \includegraphics[width=0.85\linewidth]{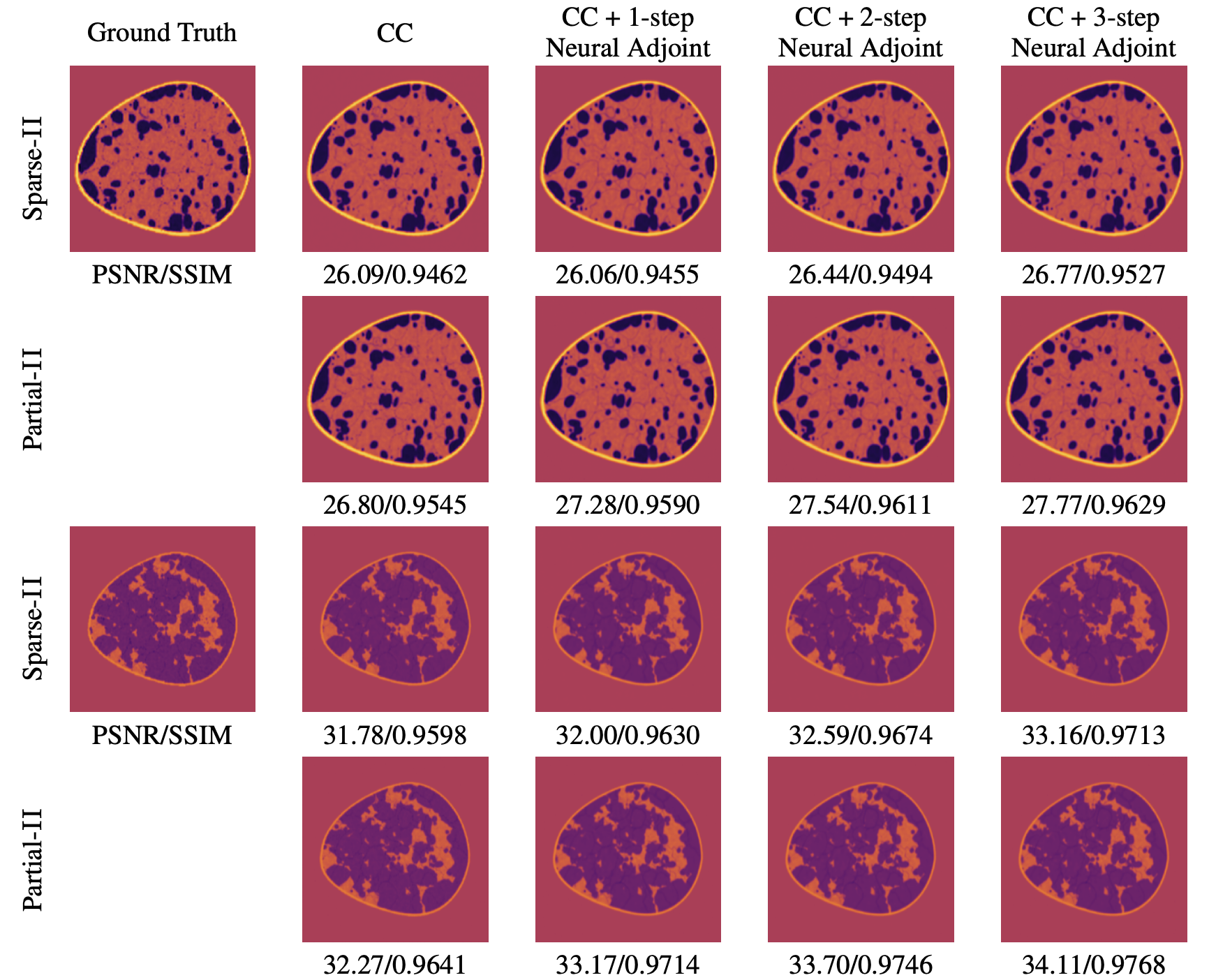}
    \caption{Visual comparison of reconstruction improvements with increasing neural adjoint-based optimization steps under severe ill-posed scenarios: Sparse-$\mathrm{II}$ and Partial-$\mathrm{II}$. The quantitative metrics PSNR/SSIM are shown below.}
   \label{fig:multistep}
\end{figure}

In \cref{fig:multistep}, it visually demonstrates that each additional optimization step enhances anatomical details including subtle heterogeneities and internal textures. Although most samples show strictly progressive improvement, we observe that some samples show a slight degradation ($0.2$-$0.5$ dB PSNR decrease) at the first sampling step $\tau_5$. This initial drop likely results from the neural adjoint’s guidance disrupting the sampling trajectory of the conditional CM at the initial sampling phase. The quantitative results in \cref{tab:multistep} show that the performance improves monotonically with more optimization steps, across all noise levels. 
\begin{table}[htp!]
\renewcommand{\arraystretch}{1.5}  
\centering
\resizebox{\textwidth}{!}{%
\begin{tabular}{l@{\hskip 10pt}ccc@{\hskip 20pt}ccc}
\hline
& \multicolumn{3}{c}{\textbf{Sparse-$\mathrm{II}$}} & \multicolumn{3}{c}{\textbf{Partial-$\mathrm{II}$}} \\
\cline{2-7}
\textbf{Methods} & \textbf{Noise-Free} & \textbf{Noise-10dB} & \textbf{Noise-5dB}& \textbf{Noise-Free} & \textbf{Noise-10dB} & \textbf{Noise-5dB}  \\ 
\hline
Conditional CM & 28.13 / 0.9279 & 27.51 / 0.9169 & 26.18 / 0.8901 & 28.15 / 0.9266 & 27.83 / 0.9200 & 27.01 / 0.9045  \\
Conditional CM + Neural Adjoint Optimization 1-step & 28.84 / 0.9378 & 28.27 / 0.9291 & 26.84 / 0.9044 & 28.61 / 0.9337 & 28.33 / 0.9285 & 27.62 / 0.9159  \\
Conditional CM + Neural Adjoint Optimization 2-step & \underline{29.31} / \underline{0.9434} & \underline{28.75} / \underline{0.9358} & \underline{27.38} / \underline{0.9140} & \underline{28.85} / \underline{0.9366} & \underline{28.56} / \underline{0.9314} & \underline{27.94} / \underline{0.9206}  \\
Conditional CM + Neural Adjoint Optimization 3-step& \textbf{29.75} / \textbf{0.9486} & \textbf{29.22} / \textbf{0.9421} & \textbf{27.91} / \textbf{0.9231} & \textbf{29.02} / \textbf{0.9392} & \textbf{28.71} / \textbf{0.9346} & \textbf{28.06} / \textbf{0.9232}  \\
\hline
\end{tabular}
}
\caption{Quantitative evaluation of neural adjoint optimization steps across different noise levels. The table reports average PSNR/SSIM values, demonstrating consistent improvement with more optimization steps. \textbf{Bold} and \underline{underline} entries indicate best and second-best performances respectively. }
\label{tab:multistep}
\end{table} 

\subsection{Quantitative Results}\label{sec:6.3}

\begin{table}[htp!]
\renewcommand{\arraystretch}{1.5}  
\centering
\resizebox{\textwidth}{!}{%
\begin{tabular}{l@{\hskip 10pt}cccc@{\hskip 20pt}cccc}
\hline
& \multicolumn{4}{c}{\textbf{Sparse-View}} & \multicolumn{4}{c}{\textbf{Partial-View}} \\
\cline{2-9}
& \multicolumn{2}{c}{\textbf{Config  $\mathrm{I}$}}  & \multicolumn{2}{c}{\textbf{Config $\mathrm{II}$}} &\multicolumn{2}{c}{\textbf{Config  $\mathrm{I}$}} &\multicolumn{2}{c}{\textbf{Config  $\mathrm{II}$}}\\
\cline{2-9}
\textbf{Methods} & \textbf{PSNR $\uparrow$} & \textbf{SSIM $\uparrow$} & \textbf{PSNR $\uparrow$} & \textbf{SSIM $\uparrow$} & \textbf{PSNR $\uparrow$} & \textbf{SSIM $\uparrow$} & \textbf{PSNR $\uparrow$} & \textbf{SSIM $\uparrow$} \\ 
\hline

\multicolumn{7}{l}{\textbf{Traditional Method}} \\
\hline
Convergent Born Series (CBS-solver) \cite{osnabrugge2016convergent} & 21.05 & 0.5891 & 18.26 & 0.2054  & 20.41 & 0.5080 & 17.74 & 0.1332 \\

\hline
\multicolumn{7}{l}{\textbf{Unsupervised Diffusion-based Sampling}} \\
\hline
Diffusion Posterior Sampling (DPS) \cite{chung2022diffusion} & 24.36 & 0.7873 & 22.16 & 0.6652  & 21.09 & 0.5921 & 19.36 & 0.5032 \\
Decomposed Diffusion Sampler (DDS) \cite{chung2023decomposed} & 27.98 & 0.9028 & 25.80 & 0.8313  & 25.48 & 0.8425 & 21.95 & 0.6348 \\
\hline   
\multicolumn{7}{l}{\textbf{Supervised End-to-End Networks }} \\
\hline
FNO-based Inversion (NIO) \cite{molinaro2023neural} & 26.30 & 0.8904 & 23.33 & 0.7798  & 26.39 & 0.8931 & 24.89 & 0.8362 \\
CNN-based Inversion (Inversion-Net) \cite{zeng2021inversionnet3d} & 28.21 & 0.9335 & 27.56 & 0.9208 & 27.15 & 0.9162 & 26.17 & 0.8960 \\
\hline
\multicolumn{7}{l}{\textbf{Diff-ANO}} \\
\hline
Conditional Tweedie + Adjoint Neural Operator & \underline{31.37} & \underline{0.9678} & \underline{29.72} & \underline{0.9477} & \underline{30.25} & \underline{0.9584} & \underline{28.44} & \underline{0.9329} \\
Conditional CM + Adjoint Neural Operator & \textbf{32.07} & \textbf{0.9732} & \textbf{30.24} & \textbf{0.9576} & \textbf{30.42} & \textbf{0.9615} & \textbf{29.51} & \textbf{0.9510} \\
\hline
\end{tabular}}
\caption{Quantitative comparison (PSNR/SSIM) of reconstruction methods under sparse-view and partial-view measurement configurations with three noise levels. \textbf{Bold} and \underline{underline} denote best and second-best results respectively.
}
\label{tab:quant_results}
\end{table}

In this section, we present a comparison of quantitative results of various methods under different measurement scenarios, including sparse-view and partial-view settings with three noise levels. Our proposed framework is evaluated in two distinct configurations, both employing a $5$-step neural adjoint optimization scheme as described in \cref{sec:6.1}. The primary difference between them lies in their consistency sampling strategies: 

\begin{enumerate}
     \item \textit{Conditional Tweedie + Adjoint Neural Operator}: employs fixed step-size gradient descent ($\eta = 0.1$) for optimization, combined with the conditional Tweedie model $
     \mathbb{E}_{\Theta^{*}_{\rm CC}}[\x_{0}|\x_{t}, \y^\delta]$ as the multi-step sampling strategy.
     \item \textit{Conditional CM + Adjoint Neural Operator}: utilizes the same optimization parameters but implements the conditional consistency model $f_{\Theta^{*}_{\rm CC}}(\x_{t}, \y^\delta, t)$ for multi-step sampling.
\end{enumerate}
In \cref{tab:quant_results}, the proposed configurations demonstrate superior reconstruction fidelity compared to the traditional method, unsupervised diffusion-based methods, and supervised end-to-end networks. The CBS-Solver exhibits fundamental limitations in handling measurement-induced ill-posedness without explicit prior regularization, resulting in degraded reconstruction quality across all scenarios. Unsupervised diffusion methods (DPS and DDS) show marked improvement through diffusion priors, achieving reasonable performance on testing samples.  Although direct learning-based approaches (NIO and Inversion-Net) achieve competitive results in all measurement settings, their performance suffers from generalization limitations inherent to the data-driven frameworks. Notably, even our method variant $\mathrm{I}$ outperforms all comparative methods through its integration of adjoint neural operators with unconditional consistency priors. Variant $\mathrm{II}$ achieves additional performance gains by incorporating conditional consistency sampling, demonstrating: (1) the efficacy of data-driven conditioning in consistency models, and (2) the benefit of combining adjoint-based optimization with learned neural operators.

\subsection{Computational Efficiency}\label{sec:6.4}

In \cref{tab:evaluation_efficiency}, we evaluate the computational efficiency of our proposed method against baseline algorithms under sparse-$\mathrm{I}$ scenario. Specifically, we focus on two primary metrics to assess computational cost: (1) the number of neural network evaluations (NFE), and (2) the number of PDE evaluations (NPE). The NFE indicates how many times neural networks—including score-based and consistency-based models—are evaluated, whereas NPE represents the number of PDE solves, mainly associated with neural operators and CBS/BCBS-solvers. Notably, each execution of the neural operator or CBS/BCBS-solver for predicting Helmholtz wavefields across all sources counts as $1$ NPE. Since we utilize \cref{eq:linear_adjoint} to avoid explicitly solving the adjoint Helmholtz equation, thus only single evaluation is required for each adjoint-based gradient calculation.

\begin{table}[htp!]
\renewcommand{\arraystretch}{1.4}
\centering
\scalebox{0.9}{
\begin{tabular}{c|cccccccc}
\hline
& CBS-solver \cite{osnabrugge2016convergent} (CPU)& CBS-solver (GPU)& DPS \cite{chung2022diffusion} & DDS \cite{chung2023decomposed}  & Diff-ANO\\
\hline
NFE  & 0 & 0 & 1000 & 50  & 5 \\
NPE  & 60+ & 60+ & 1000 & 100+ & 5 \\
Average Time  & 10.9h & 298.5s & 4023.7s & 451.9s & 1.1s  \\
\hline
\end{tabular}
}

\caption{Computational efficiency comparison of reconstruction methods: average time per sample under sparse-$\mathrm{I}$ scenario. NFE denotes the number of neural network/operator evaluations, and NPE represents PDE evaluations required by each method.}
\label{tab:evaluation_efficiency}
\end{table}

The CBS-solver, executed using iterative loops for the Helmholtz equation, employs L-BFGS with $30$ optimization steps. Given that L-BFGS requires additional forward evaluations for line search, the NPE surpasses $60$ evaluations per reconstruction. Conversely, the BCBS-solver implements batching on A100 GPUs to accelerate PDE solves, dramatically reducing computational time from 10.9 hours (CBS-solver) to 298.5 seconds. This highlights that while Diff-ANO runs on GPUs, the baseline CBS is heavily bottlenecked by CPU execution, strengthening the need for neural acceleration.

It is important to address the training costs, which are substantial but ultimately amortized. The pipeline requires EDM pretraining, consistency distillation, ControlNet fine-tuning, and crucially, MgNO training (which alone costs approximately 130 A100 GPU-hours due to the on-the-fly BCBS data generation). While the upfront computational investment is high, these models scale efficiently during inference. In clinical or continuous monitoring scenarios where thousands of reconstructions are needed for a fixed measurement geometry, the $1.1$s inference time justifies the initial training cost. The DPS method, despite utilizing the GPU-accelerated BCBS-solver, employs $1000$ discretized sampling steps coupled with traditional gradient descent, resulting in an equal count with both NFE and NPE. DDS partially alleviates this computational burden by integrating $50$-step DDIM accelerated sampling and conjugate gradient descent. Nevertheless, due to the iterative nature of the conjugate gradient algorithm necessitating line searches, DDS still incurs over $100$ PDE evaluations. Our proposed approach circumvents direct PDE solving by transforming all PDE-related computations into neural network evaluations, substantially increasing computational efficiency.  Within the conditional CM framework, our method only requires a few-step evaluations ($5$ neural operator evaluations and $5$ consistency model evaluations) to perform measurement-constrained iterative refinement. Consequently, our method achieves an impressive acceleration, reducing the computational time per sample to merely $1.1$ seconds—orders of magnitude faster than all compared methods—while maintaining high-quality reconstruction.

\subsection{Multi-frequency Reconstruction Results}\label{sec:6.5}

In traditional USCT, multi-frequency reconstruction (also known as frequency hopping or frequency continuation) is a pivotal strategy to mitigate the strong nonlinearity and cycle-skipping artifacts inherent in wave-based inverse problems. Lower frequencies provide a smoother, more convex objective function that captures large-scale background structures, while higher frequencies incrementally introduce fine-scale resolution at the cost of increased non-convexity.

To evaluate the compatibility of Diff-ANO with this classical acoustic protocol, we conducted experiments across a frequency spectrum from \SI{400}{kHz} to \SI{600}{kHz} with a \SI{50}{kHz} step size. At each frequency stage, we execute 5 sequential neural adjoint refinements, transitioning through $\tau = \{0.1, 0.12, 0.14, 0.16, 0.18\}$ to progressively align the reconstruction with the measurement-consistent manifold. The quantitative performance under the highly ill-posed sparse-view $\mathrm{II}$ (\SI{5}{dB} SNR) is presented in \cref{tab:multi_freq}.

\begin{table}[htp!]
\renewcommand{\arraystretch}{1.2}
\centering
\scalebox{0.9}{
\begin{tabular}{lcccc}
\toprule
&\multicolumn{2}{c}{Single \SI{500}{kHz} }& \multicolumn{2}{c}{\SI{400}{kHz} $\sim$ \SI{600}{kHz}}\\\cline{2-5}
\textbf{Optimization Strategy} & \textbf{PSNR $\uparrow$} & \textbf{SSIM $\uparrow$}& \textbf{PSNR $\uparrow$} & \textbf{SSIM $\uparrow$} \\\midrule
L-BFGS & 18.26  & 0.2054& 23.11 & 0.7230 \\
Diff-ANO & 30.24 & 0.9576 & \textbf{31.84} & \textbf{0.9719}\\
\bottomrule
\end{tabular}
}
\caption{Quantitative comparison of reconstruction strategies under Sparse-view $\mathrm{II}$ (\SI{5}{dB} SNR) using single and multi-frequency schedules.}
\label{tab:multi_freq}
\end{table}
\paragraph{Analysis and Discussion.} 
The results in \cref{tab:multi_freq} demonstrate a profound synergy between multi-scale physics and generative regularization. For the traditional L-BFGS method, the transition from single-frequency to multi-frequency optimization leads to a significant performance leap (from \SI{18.26}{dB} to \SI{23.11}{dB} PSNR). This gain confirms that the sequential inclusion of low-to-high frequency data helps the gradient-based solver escape poor local minima by "warming up" the acoustic background. However, even with multi-frequency information, L-BFGS remains limited by the sparse-view configuration, which manifests as persistent oscillatory artifacts and low SSIM.

Crucially, when equipped with the multi-frequency schedule, Diff-ANO reaches a peak performance of \SI{31.84} dB PSNR and \SI{0.9719} SSIM. This improvement can be attributed to two factors: 
(1) \textit{Spectral Redundancy as Regularization}: The multi-frequency approach provides a wider spectral coverage of the scattering field, which MgNO successfully exploits to stabilize the gradient updates across different scales. 
(2) \textit{Progressive Manifold Alignment}: The 5-step refinement per frequency ($\tau = 0.2 \to 0.1$) allows the model to reconcile the high-level data-driven prior with the frequency-specific wave physics. This demonstrates that Diff-ANO does not merely "memorize" the prior but dynamically adapts it to conform to multi-scale physical observations.

\subsection{Results on OpenPros Dataset}\label{sec:6.6}

To further evaluate the generalization capability and robustness of Diff-ANO in clinically realistic scenarios, we extend our experiments to the OpenPros  dataset \cite{wang2025openpros}. OpenPros is the first large-scale benchmark specifically designed for limited-angle prostate USCT, featuring anatomically precise 2D speed-of-sound phantoms derived from clinical MRI/CT scans and ex vivo specimens. These phantoms incorporate high tissue heterogeneity and the presence of pelvic bones, which present significant challenges for wave-based reconstruction.

\begin{table}[htp!]
\renewcommand{\arraystretch}{1.2}
\centering
\scalebox{0.9}{
\begin{tabular}{lcccc}
\toprule
\textbf{Method} &L-BFGS &DDS  &Inversion-Net &Diff-ANO\\ \midrule
\textbf{PSNR $\uparrow$} &18.54 &25.12 &28.35 &\textbf{30.47}\\ \midrule
\textbf{SSIM $\uparrow$}      &0.5213& 0.8145 & 0.9212 & \textbf{0.9588}\\ \bottomrule
\end{tabular}
}
\caption{Quantitative comparison of reconstruction methods on the OpenPros subset (\SI{10}{dB} SNR) under the limited-angle dual-probe configuration.}
\label{tab:openpros_results}
\end{table}

In this study, we randomly selected a subset of \SI{8000} 2D-SOS images from the OpenPros dataset. Rather than using the original time-domain full-waveform data provided in the dataset, we utilized the high-fidelity SOS maps to simulate measurements governed by the Helmholtz equation at \SI{500}{kHz} using our CBS-based forward solver. To emulate realistic clinical conditions, we introduced \SI{10}{dB} additive Gaussian noise to the measurements. The acquisition geometry strictly follows the dual-probe configuration described in \cite{wang2025openpros}, which inherently creates a severely limited-angle observation window. The quantitative results comparing traditional L-BFGS, unsupervised DDS, supervised Inversion-Net, and our Diff-ANO are summarized in \cref{tab:openpros_results}.

\textbf{Remarks.} The CBS-solver (CPU) for multi-source simulations utilizes \texttt{Intel Xeon Platinum 8358P} CPUs, while the CBS-solver (GPU) leverage \texttt{NVIDIA A100} GPUs for batched Helmholtz solves. Besides, the training and inference processes of both conditional consistency models and neural operators are implemented in PyTorch 1.13, trained/evaluated on $4\times$\texttt{A100} GPUs with parallelism.

\section{Ablation Study}\label{sec:7}

\subsection{Inversion Blocks for Conditioning Consistency Model}\label{sec:7.1}
Our ablation study evaluates three inversion blocks for initial estimation in the conditional CM. In \cref{tab:inversion_blocks}, one-step GD leverages the precomputed background Helmholtz solutions without real-time PDE solving. The trained direct-mapping networks (NIO/Inversion-Net) are directly plugged in as inversion blocks to evaluate the performance.  Inversion-Net delivers more optimal results than NIO, aligning with our algorithm’s emphasis on high-fidelity pre-reconstruction in \cref{sec:4.2}, where better initial estimates enable more effective refinement through measurement-conditioned consistency constraints.

\begin{table}[htp!]
\renewcommand{\arraystretch}{1.5}  
\centering
\scalebox{0.9}{
\begin{tabular}{l@{\hskip 20pt}cc@{\hskip 20pt}cc}
\hline
& \multicolumn{2}{c}{{\hskip -10pt}\textbf{Sparse-$\mathrm{II}$}} & \multicolumn{2}{c}{\textbf{Partial-$\mathrm{II}$}} \\
\cline{2-5}
\textbf{Inversion Blocks} & PSNR$\uparrow$ & SSIM$\uparrow$ & PSNR$\uparrow$ & SSIM$\uparrow$  \\ 
\hline
One-step GD & 23.76 & 0.7956 & 22.28 &  0.7766  \\
NIO \cite{molinaro2023neural} & 27.48 & 0.9187 & 28.29 & 0.9206  \\
Inversion-Net \cite{zeng2021inversionnet3d} & 30.24 & 0.9576  & 29.51 & 0.9510   \\
\hline
\end{tabular}
}
\caption{Comparative evaluation of inversion blocks in sparse-$\mathrm{II}$ and partial-$\mathrm{II}$ scenarios. The quantitative results (PSNR, SSIM) of our framework are evaluated for different inversion blocks with the fixed neural operator (MgNO). }
\label{tab:inversion_blocks}
\end{table}

\subsection{Neural Operators for Adjoint-Based Optimization}\label{sec:7.2}
\begin{figure}
    \centering
    \includegraphics[width=0.875\linewidth]{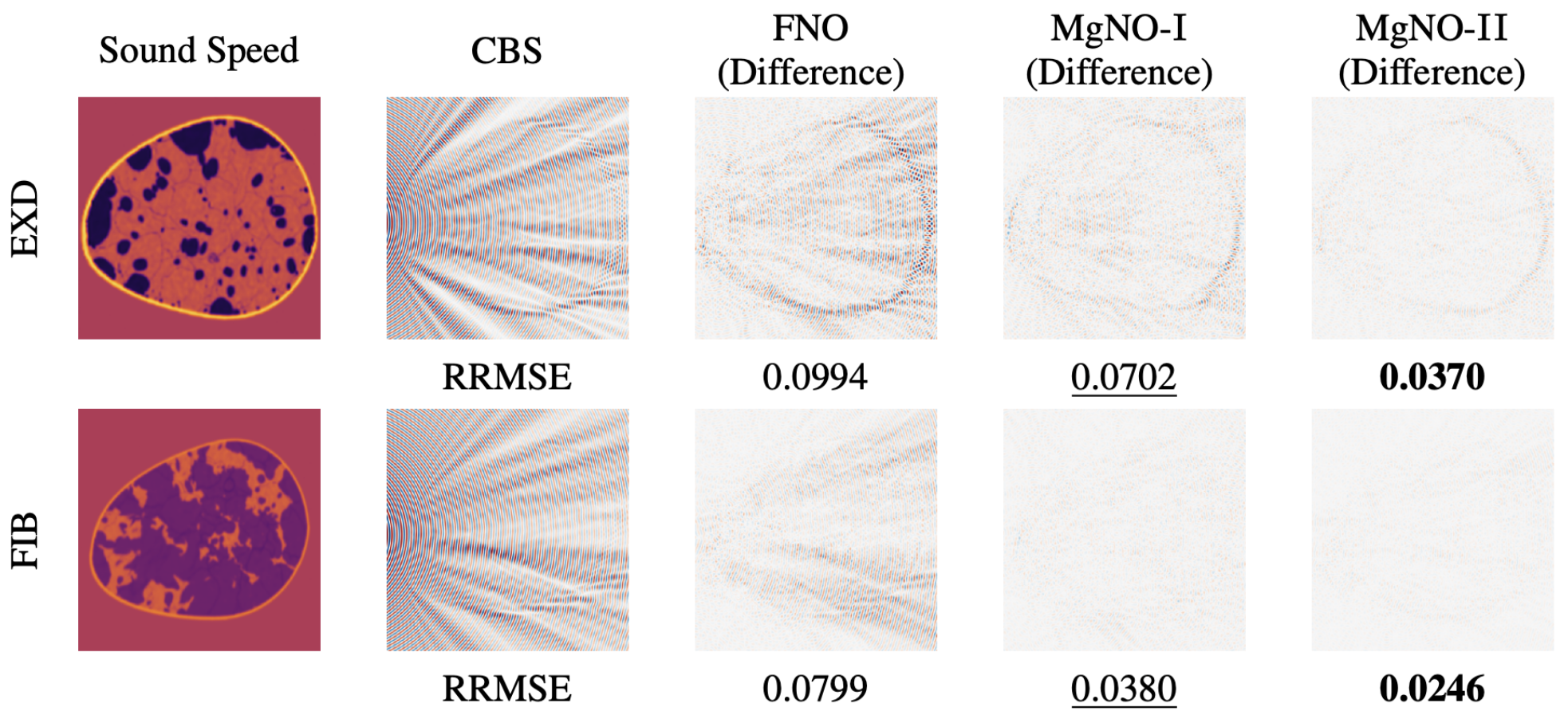}
    \caption{Visual comparison of forward wavefield prediction differences for EXD-type and FIB-type samples. Sound-speed fields and the corresponding CBS-solver wavefield results are shown as references. }
\label{fig:operator_ablation}
\end{figure}

In \cref{fig:operator_ablation} and \cref{tab:evaluation}, it reveals the impact of neural operator architecture on Helmholtz-based forward prediction and inversion. Here, FNO \cite{li2020fourier} serves as the baseline neural operator for approximating the forward Helmholtz operator. In our implementations, FNO utilizes $4$ spectral convolution layers, configured with $25$ Fourier modes and a channel width of $32$. MgNO-$\mathrm{I}$ uses a smaller feature dimension ($12$ channels for multi-scale layers) and fewer recurrent iterations ($4$ applications of the V-cycle), making it a lightweight version of MgNO. MgNO-$\mathrm{II}$, on the other hand, utilizes a larger feature dimension ($24$ channels) and a higher number of recurrent iterations ($6$ applications of the V-cycle). The increased channels and iterations of MgNO-$\mathrm{II}$ allow it to better capture the solution's characteristics, resulting in superior performance over both FNO and MgNO-$\mathrm{I}$ for forward prediction and inversion. Notably, although the neural operator's accuracy determines the adjoint-based gradient reliability, even the inaccurate neural operator (FNO) can enhance the performance via the adjoint-based optimization, as shown in \cref{tab:evaluation}. 

\begin{table}[htp!]
\renewcommand{\arraystretch}{1.5}  
\centering
\scalebox{0.9}{
\begin{tabular}{lcccccc}
\hline
&\textbf{Neural Operators}  &FNO&MgNO-$\mathrm{I}$ &MgNO-$\mathrm{II}$ & CBS-solver \\\hline
\textbf{Forward Prediction}&RRMSE&0.0882 &0.0413&0.0264&--\\\hline
\multirow{2}{*}{\textbf{Inversion}}&PSNR&29.58&30.96&32.07&32.85\\
&SSIM&0.9520&0.9633&0.9732&0.9803\\
\hline
\end{tabular}
}
\caption{The Relative Root Mean Square Error (RRMSE) of forward prediction is evaluated for different neural operators. Under sparse-$\mathrm{I}$ scenario, PSNR/SSIM metrics of inversion are evaluated for different neural operators and baseline CBS-solver with the fixed inversion block (Inversion-Net).} 
\label{tab:evaluation}
\end{table}

\subsection{Limitations}\label{sec:8.2}

\paragraph{Supervised Guidance for Conditional CMs.} The proposed method relies on an appropriate initial reconstruction with the iterative PDE‐based refinement, owing to the ill‐posedness introduced by the under-sampled measurement configuration. In our framework, the conditional CM furnishes this initial estimate and, in the multi‐step sampling scheme, seamlessly integrates physics‐informed guidance via the adjoint neural operator. However, unlike unsupervised sampling strategies (e.g., DPS or unconditional CM), our paradigm necessitates paired datasets of under-sampled measurements and ground‐truths for conditional CM training. While training the neural operator does not depend on a specific measurement setup, the overall pipeline remains a supervised-learning method constrained by the configuration‐specific training data.

\paragraph{Dependence on Self-Adjoint Structure.} A key enabler of the adjoint neural operator is the \textit{self‐adjointness} of the Helmholtz operator, allowing the adjoint solution \cref{eq:adjoint} to be expressed as a linear combination of forward Helmholtz solutions \cref{eq:linear_adjoint}. Consequently, training surrogate forward neural operators suffices for efficient gradient optimization. However, extending the proposed framework to a broader research direction exploring complex non-self-adjoint transport and wave phenomena requires careful theoretical analysis to independently map the forward and adjoint operators, necessitating co-designed network architectures and distinct training strategies.

\section{Conclusion}\label{sec:8}
In this work, we have presented a novel hybrid framework for USCT by integrating a conditional consistency model with adjoint neural operators. Our approach departs from the conventional adjoint-based methods that rely heavily on numerical PDE solvers, and instead leverages neural operators for gradient computation in the consistency sampling framework. Through extensive experiments, we demonstrate that our method achieves high-fidelity reconstructions in only a few sampling steps, significantly accelerating USCT imaging. By bridging the continuous domain of PDE-constrained optimization with the discrete framework of generative priors, we construct a unified pipeline that fundamentally overcomes the computational bottlenecks of traditional solvers. This integration provides a powerful paradigm that combines data-driven priors with physics-based constraints for efficient and high-fidelity USCT reconstruction.
\bibliographystyle{unsrt}  
\bibliography{references}

\end{document}